\documentclass[11pt,letterpaper]{amsart}
\usepackage{amsmath,amssymb,color,tikz-cd}
\usepackage[backref,pagebackref,pdftex,hyperindex]{hyperref}
\usepackage[alphabetic]{amsrefs}
\usepackage{amssymb}
\usepackage{bbm}
\usepackage{enumitem}
\usepackage{upgreek}

\newtheorem{proposition}{Proposition}[section]
\newtheorem{theorem}[proposition]{Theorem}
\newtheorem{lemma}[proposition]{Lemma}
\newtheorem{corollary}[proposition]{Corollary}

\theoremstyle{definition}
\newtheorem{remark}[proposition]{Remark}
\newtheorem{definition}[proposition]{Definition}

\title{On the shape of the K-semistable domain and wall crossing for K-stability}
\author{Chuyu Zhou}
\address{School of Mathematical Sciences, Xiamen University, Siming South Road 422, Xiamen, Fujian 361005, China}
\email{chuyuzhou1@gmail.com}

\date{} 

\thanks{2010 
	    \emph{Mathematics Subject Classification}: 14J45.
	    \newline
	    \indent 
		\emph{Keywords}: Log Fano pair, K-stability, K-semistable domain, K-moduli, wall crossing.
        \newline
		\indent
		\emph{Competing interests}: The author declares none.
		}

\newcommand{\ord}{{\rm {ord}}}
\newcommand{\tc}{{\rm {tc}}}
\newcommand{\vol}{{\rm {vol}}}

\newcommand{\red}{{\rm {red}}}

\newcommand{\Val}{{\rm {Val}}}

\newcommand{\Kss}{{\rm {Kss}}}
\newcommand{\Cone}{{\rm {Cone}}}
\newcommand{\LC}{{\rm {LC}}}

\newcommand{\bA}{\mathbb{A}}

\newcommand{\bC}{\mathbb{C}}

\newcommand{\bN}{\mathbb{N}}

\newcommand{\bP}{\mathbb{P}}
\newcommand{\bQ}{\mathbb{Q}}
\newcommand{\bR}{\mathbb{R}}

\newcommand{\mB}{\mathcal{B}}

\newcommand{\mD}{\mathcal{D}}
\newcommand{\mE}{\mathcal{E}}

\newcommand{\mG}{\mathcal{G}}
\newcommand{\mH}{\mathcal{H}}

\newcommand{\mL}{\mathcal{L}}
\newcommand{\mM}{\mathcal{M}}

\newcommand{\mO}{\mathcal{O}}
\newcommand{\mP}{\mathcal{P}}

\newcommand{\mX}{\mathcal{X}}
\newcommand{\mY}{\mathcal{Y}}

\begin{document}

\begin{abstract}
Fixing two positive integers $d$ and $k$, a positive number $v$, and a positive integer $I$, we prove that the K-semistable domain of the log pair $(X, \sum_{j=1}^kD_j)$ is a rational polytope lying in the $k$-dimensional simplex $\overline{\Delta^k}$, where $X$ is a Fano variety of dimension $d$, $D_j\sim_\mathbb{Q} -K_X$, $(-K_X)^d=v$, $I(K_X+D_j)\sim 0$, and $(X, \sum_{j=1}^kc_jD_j)$ is a K-semistable log Fano pair for some $c_j\in [0,1)\cap \mathbb{Q}$. Moreover, we show that there are only finitely many polytopes  which may appear as the K-semistable domains for such log pairs. Based on this, we establish a wall crossing theory for K-moduli with multiple boundaries.

\end{abstract}

\maketitle

\setcounter{tocdepth}{1}

\tableofcontents

\section{Introduction}

We work over the complex number field $\bC$ throughout the article.

In the past few years, a rather complete theory on algebraic K-stability has been developed, which finally leads to the construction of a projective moduli space for Fano varieties with K-stability (e.g. \cite{Xu21}). One novel characteristic of K-stability is the wall crossing phenomenon (e.g. \cite{ADL19, Zhou23}), which provides the suitable ideas and tools to study the birational geometry of various kinds of moduli spaces (e.g. \cite{ADL19, ADL20, ADL21}). To study the wall crossing phenomenon for K-stability, an important ingredient is to confirm the K-semistable domain. We start with the following definition of K-semistable domain in an easier case.

\begin{definition}
Let $(X, D)$ be a log pair, where $X$ is a Fano variety of dimension $d$ and $0\leq D\sim_\bQ -K_X$. The \textit{K-semistable domain} of $(X, D)$ is defined as follows:
$$\Kss(X, D):=\overline{\{t\in [0, 1)\cap \bQ\ |\ \textit{$(X, tD)$ is K-semistable}\}}, $$
where the overline means taking the closure.
\end{definition}

Let $(X, D)$ be a log pair as in the above definition, we define
$$l(X, D):=\inf \{t\in [0,1)\cap \bQ\ |\ \textit{$(X, tD)$ is K-semistable}\}, $$
$$u(X, D):=\sup \{t\in [0,1)\cap \bQ\ |\ \textit{$(X, tD)$ is K-semistable}\}, $$
then it is clear that $\Kss(X, D)=[l(X, D), u(X, D)]$. Actually, we have the following more general result:

\begin{theorem}{\rm{(\cite{Zhou23})}}\label{thm: 1 component}
Fix a positive integer $d$, a positive number $v$, and a finite set $\Phi$ of non-negative rational numbers. Then there exist finite rational numbers
$$0=c_0<c_1<...<c_k<c_{k+1}=1 $$
depending only on $d, v, \Phi$ such that for any log pair $(X, D)$ satisfying
\begin{enumerate}
\item $X$ is a Fano variety of dimension $d$ and $(-K_X)^d=v$,
\item $0\leq D\sim_\bQ -K_X$ and the coefficients of $D$ are contained in $\Phi$,
\item $(X, tD)$ is K-semistable for some $t\in [0,1)\cap \bQ$,
\end{enumerate}
we have $\Kss(X, D)=[c_i, c_j]$ for some $0\leq i\leq j\leq k+1$.
\end{theorem}

Based on this finite chamber decomposition of $[0,1]$, one could naturally establish a wall crossing theory for K-moduli parametrizing the log pairs $(X, cD)$ via changing $c$, where $(X, D)$ satisfies the conditions in Theorem \ref{thm: 1 component} and $c\in [0,1)\cap \bQ$ (e.g. \cite{ADL19, Zhou23}).
The main result in this article is to generalize Theorem \ref{thm: 1 component} to the case where we allow multiple boundaries, which is conjectured in \cite[Conjecture 8.3]{Zhou24}. Before stating the theorem, let us first recall the following set of log pairs introduced in \cite{LZ23}.

\begin{definition}
Fixing two positive integers $d$ and $k$, a positive number $v$,  and a positive integer $I$, we consider the set $\mE:=\mE(d,k,v, I)$ of log pairs $(X, \sum_{i=j}^k D_j)$ satisfying the following conditions:
\begin{enumerate}
\item $X$ is a Fano variety of dimension $d$ and $(-K_X)^d=v$;
\item $D_j$ is an effective $\bQ$-divisor satisfying $D_j\sim_\bQ -K_X$ for every $j$;
\item $I(K_X+D_j)\sim 0$ for every $j$;
\item there exists a vector $(c_1,...,c_k)\in \Delta^k$ such that $(X, \sum_{j=1}^k c_jD_j)$ is K-semistable, where 
$$\Delta^k:=\{(c_1,...,c_k)\ |\ \textit{$c_j\in [0,1)\cap \bQ$ and $0\leq \sum_{j=1}^k c_j<1$}\}.$$
\end{enumerate}
We define the \textit{K-semistable domain} of $(X, \sum_{j=1}^k D_j)\in \mE$ as follows:
$$\Kss(X, \sum_{j=1}^k D_j):=\overline{\{(a_1,...,a_k)\in \Delta^k\ |\ \textit{$(X, \sum_{j=1}^k a_jD_j)$ is K-semistable}\}}, $$
where the overline means taking the closure.
\end{definition}

We have the following characterization of K-semistable domains of log pairs in $\mE$.

\begin{theorem}\label{thm: main1}
Let $(X, \sum_{j=1}^k D_k)$ be a log pair in $\mE:=\mE(d,k, v, I)$, then
\begin{enumerate}
\item $\Kss(X, \sum_{j=1}^kD_j)$ is a polytope in the simplex $\overline{\Delta^k}$;
\item all the vertices of $\Kss(X, \sum_{j=1}^kD_j)$ are rational, which means that their coordinates are rational numbers;
\item there are only finitely many polytopes in $\overline{\Delta^{k}}$ which could appear as the  K-semistable domains of log pairs in $\mE$.
\end{enumerate}
\end{theorem}

Throughout the paper, a polytope is always assumed to be convex and compact. Note that in the statement of Theorem \ref{thm: main1}, the dimension of the K-semistable domain is not necessarily maximal as a polytope in $\overline{\Delta^{k}}$. Suppose $P$ is a given polytope,  we use $P^\circ$ to denote the set of interior points of $P$ (under relative topology), and $P^\circ(\bQ)$ the set of rational points of $P^\circ$. Theorem \ref{thm: main1} leads to the following finite chamber decomposition of $\overline{\Delta^k}$ for log pairs in $\mE$.

\begin{corollary}\label{cor: chamber decomposition}
Fix two positive integers $d$ and $k$, a positive number $v$, and a positive integer $I$.
There exists a finite chamber decomposition $\overline{\Delta^k}=\cup_{i=1}^r P_i$, where $P_i$'s are rational polytopes and $P_i^\circ\cap P_j^\circ=\emptyset$ for $i\ne j$ , such that for any log  pair $(X, \sum_{j=1}^kD_j)\in \mE(d,k,v,I)$ and any face $F$ of any chamber $P_i$, the K-semistability of $(X, \sum_{j}^kx_jD_j)$ does not change as $\vec{x}:=(x_1,...,x_k)$ varies in $F^\circ(\bQ)$. 
\end{corollary}

The variation of K-semistability under the variation of coefficients in the above corollary can be compared to the variation of GIT-semistability under the variation of linearizations (e.g. \cite{Tha96, DH98}).

The finite chamber decomposition also implies a wall crossing theory for K-moduli. For any rational vector $\vec{c}=(c_1,...,c_k)\in \Delta^k$, we denote $\mM^K_{d,k,v,I,\vec{c}}$ to be the K-moduli stack parametrizing the following subset of $\mE$:
$$\mE_{\vec{c}}:=\{(X, \sum_{j=1}^kD_j)\in \mE|\textit{$(X, \sum_{j=1}^kc_jD_j)$ is a K-semistable log Fano pair}\}. $$
Let $M^K_{d,k,v,I,\vec{c}}$ be the corresponding K-moduli space. We have the following result (refer to Section \ref{sec: K-moduli} for detailed explanation).

\begin{theorem}\label{thm: main2}
Fix two positive integers $d$ and $k$, a positive number $v$, and a positive integer $I$. There exists a finite chamber decomposition $\overline{\Delta^k}=\cup_{i=1}^r P_i$, where $P_i$'s are rational polytopes and $P_i^\circ\cap P_j^\circ=\emptyset$ for $i\ne j$, such that  
\begin{enumerate}
\item For any face $F$ of any chamber $P_i$, $\mM^K_{d,k,v,I,\vec{c}}$ {\rm{(resp. $M^K_{d,k,v,I,\vec{c}}$)}} does not change as $\vec{c}$ varies in $F^\circ(\bQ)$; 
\item Let $P_{i}$ and $P_{j}$ be two different polytopes in the decomposition that share the same face $F_{ij}$. Suppose $\vec{w}_1\in P_i^\circ(\bQ)$ and $\vec{w}_2\in P_j^\circ(\bQ)$ satisify that the segment connecting them
intersects $F_{ij}$ at a point $\vec{w}$, then $\vec{w}$ is a rational point and we have the following commutative diagram for any rational numbers $0<t_1, t_2< 1$:
\begin{center}
	\begin{tikzcd}[column sep = 1em, row sep = 2em]
	 \mM_{d,k,v,I, t_1\vec{w}+(1-t_1)\vec{w}_1}^K \arrow[d,"",swap] \arrow[rr,""]&& \mM_{d,k,v,I, \vec{w}}^K\arrow[d,"",swap] &&\mM_{d,k,v,I, t_2\vec{w}+(1-t_2)\vec{w}_2}^K\arrow[d,""]\arrow[ll,""]\\
	 M_{d,k,v,I, t_1\vec{w}+(1-t_1)\vec{w}_1}^K\arrow[rr,""]&& M_{d,k,v,I, \vec{w}}^K&&M_{d,k,v,I, t_2\vec{w}+(1-t_2)\vec{w}_2}^K\arrow[ll,"", swap].
	 	 	 	\end{tikzcd}
\end{center}
Moreover, we have
\begin{enumerate}
\item $\mM_{d,k,v,I, t_l\vec{w}+(1-t_l)\vec{w}_l}^K$
{\rm{(resp. $M_{d,k,v,I, t_l\vec{w}+(1-t_l)\vec{w}_l}^K$)}} does not change as $t_l$ varies in $[0,1)\cap \bQ$, where $l=1,2$;
\item $ \mM_{d,k,v,I, t_1\vec{w}+(1-t_1)\vec{w}_1}^K\longrightarrow \mM_{d,k,v,I, \vec{w}}^K\longleftarrow \mM_{d,k,v,I, t_2\vec{w}+(1-t_2)\vec{w}_2}^K$ are open embeddings for $t_1, t_2\in [0,1)\cap \bQ$;
\item $ M_{d,k,v,I, t_1\vec{w}+(1-t_1)\vec{w}_1}^K\longrightarrow M_{d,k,v,I, \vec{w}}^K\longleftarrow M_{d,k,v,I, t_2\vec{w}+(1-t_2)\vec{w}_2}^K$ are projective morphisms for $t_1,t_2\in [0,1)\cap \bQ$.
\end{enumerate}
\end{enumerate}
\end{theorem}

We briefly present some ideas of the proof of Theorem \ref{thm: main1}. 
First of all, for a log pair $(X, \sum_{j=1}^kD_j)\in \mE$, we have the following interpolation property for K-stability: if $(X, \Delta_1)$ and $(X, \Delta_2)$ are both K-semistable log pairs (log Fano or log Calabi-Yau), where $\Delta_i$'s are proportional to $-K_X$, then $(X, t\Delta_1+(1-t)\Delta_2)$ is also K-semistable for any $t\in [0,1]\cap \bQ$. However, unlike the one boundary case as in Theorem \ref{thm: 1 component}, this interpolation property only tells us that $\Kss(X, \sum_{j=1}^kD_j)$ is convex, which is far from being a polytope. Let us take $k=2$, and assume $\Kss(X, D_1+D_2)$ is a two dimensional convex set. To show that 
$\Kss(X, D_1+D_2)$ is a polytope, at least we need to show that it is locally a polytope at every extreme point. Taking $\vec{w}$ to be an extreme point, if $\Kss(X, D_1+D_2)$ is not locally a polytope at $\vec{w}$, then one could easily find a sequence $\{\vec{w}_i\}_i$ tending to $\vec{w}$ such that $\vec{w}_i$ lies on the boundary of $\Kss(X, D_1+D_2)$ and the line connecting $\vec{w}_i$ and $\vec{w}$ separates $\Kss(X, D_1+D_2)$ for every $i$ (see Definition \ref{def: separation}). 
The insight here is to notice that $\Kss(X, D_1+D_2)$  is cut out by (infinite) hyperplanes given by beta invariants of prime divisors over $X$:
$$\beta_{X, x_1D_1+x_2D_2}(E)=A_X(E)-\sum_{j=1}^2\ord_E(D_j)\cdot x_j-(1-\sum_{j=1}^2x_j)S_X(E)=0. $$
By beta-invariant criterion (e.g. \cite{Li17, Fuj19}),  these hyperplanes cannot cross through the interior part of $\Kss(X, D_1+D_2)$. However, it turns out that each $\vec{w}_i$ produces a prime divisor $E_i$ over $X$ such that the hyperplane given by the beta invariant of $E_i$, i.e. $\beta_{X, x_1D_1+x_2D_2}(E_i)=0$,  passes through both $\vec{w}_i$ and $\vec{w}$ (see the argument in the proof of Lemma \ref{lem: polytope lemma}). This leads to a contradiction. For general $k$, by a local characterization of polytopes (see Lemma \ref{lem: local char polytope}), we could confirm that $\Kss(X, \sum_{j=1}^kD_j)$ is locally a polytope at every extreme point (see Lemma \ref{lem: polytope lemma}). Applying the similar argument on showing locally being a polytope, we could also confirm that the extreme points are in fact discrete (see Lemma \ref{lem: discrete lemma}). This forces the K-semistable domain to be a polytope. For the rationality part, the same idea in one boundary case to show the rationality of wall numbers also applies here to reveal that there are dense rational points on the boundary of $\Kss(X, \sum_{j=1}^kD_j)$ (Lemma \ref{lem: dense lemma}). This forces the extreme points to be rational. In \cite{LZ23}, we have proved the log boundedness of the set $\mE$ (Theorem \ref{thm: boundedness}), which is an important ingredient here to finally confirm the finiteness of polytopes. To achieve this, we are inspired by the idea in \cite{BLX22} where they prove the constructibility of delta invariants in a family. For a family of log pairs whose fibers are contained in $\mE$, suppose there is a fiberwise log resolution, we could construct a finite stratification of the base such that the fibers over each piece of the stratification admit the same K-semistable domain (Lemma \ref{lem: finite lemma}). This achieves the finiteness.

In Section \ref{sec: preliminaries}, we present some basic concepts that are applied in this paper. In Section \ref{sec: polytope}, we introduce convex sets and polytopes, and give a local characterization of polytopes (Lemma \ref{lem: local char polytope}). In Section \ref{sec: boundedness}, we include a complete proof of boundedness of the set $\mE$. In Section \ref{sec: rational polytope}, we show that the K-semistable domains of log pairs in $\mE$ are rational polytopes. In Section \ref{sec: finiteness}, we show that there are only finitely many polytopes that may appear as the K-semistable domains of log pairs in $\mE$. In the last section, we present the application of Theorem \ref{thm: main1} on the wall crossing theory for K-moduli.

\noindent
\subsection*{Acknowledgement}
The author is supported by the grant of European Research Council (ERC-804334).

\section{Preliminaries}\label{sec: preliminaries}

We say that $(X,\Delta)$ is a \emph{log pair} if $X$ is a normal projective variety and $\Delta$ is an effective $\bQ$-divisor on $X$ such that $K_X+\Delta$ is $\bQ$-Cartier.  The log pair $(X,\Delta)$ is called \emph{log Fano} if it admits klt singularities and $-(K_X+\Delta)$ is ample; if $\Delta=0$, we just say $X$ is a \emph{Fano variety}. 
The log pair $(X,\Delta)$ is called \emph{log Calabi-Yau} if $K_X+\Delta\sim_\bQ 0$. We will also encounter log pairs with $\bR$-coefficients, where we just replace the $\bQ$-Cartier condition with the $\bR$-Cartier condition and preserve all the other conditions. When we mention the term \textit{log pair} it always refers to the $\bQ$-coefficients case, and we will emphasize in a bracket if there are $\bR$-coefficients.
For various types of singularities in birational geometry, e.g.  klt, lc, and plt singularities, we refer to \cite{KM98,Kollar13}.

\subsection{Invariants associated to log pairs}

\subsubsection{Log discrepancy}
Let $(X,\Delta)$ be a log pair. Suppose $f\colon Y\to X$ is a proper birational morphism between normal varieties and $E$ is a prime divisor on $Y$, we say $E$ is a \textit{prime divisor over $X$}.  The following invariant associated to $E$,
$$A_{X,\Delta}(E):=1+\ord_E(K_Y-f^*(K_X+\Delta)), $$
is called the \emph{log discrepancy} of $E$ with respect to the log pair $(X,\Delta)$. Clearly one could define log discrepancies of divisors with respect to log pairs with $\bR$-coefficients. We also mention here that prime divisors over $X$ just correspond to divisorial valuations on the function field of $X$, i.e. $K(X)$, and we can define log discrepancy for any valuation $v\in \Val_X$ (e.g. \cite{JM12}), where $\Val_X$ means the valuation space over $X$.

\subsubsection{Averaging vanishing order}
Let $(X,\Delta)$ be a log Fano pair and $E$ a prime divisor over $X$. The averaging vanishing order of $E$ with respect to $(X, \Delta)$ is defined as follows:
$$S_{X,\Delta}(E):=\frac{1}{\vol(-K_X-\Delta)}\int_0^\infty \vol(-f^*(K_X+\Delta)-tE){\rm{d}}t .$$
To see why it is called averaging vanishing order, we introduce another notion. For any divisible positive integer $m$, we say $0\leq B_m\sim_\bQ -(K_X+\Delta)$ is an \textit{$m$-basis type divisor} if it is of the following form
$$B_m=\frac{\sum_{i=1}^{N_m}\{s_i=0\}}{mN_m}, $$
where $N_m=\dim H^0(X, -m(K_X+\Delta))$ and $\{s_i\}_i$ is a complete basis of $H^0(X, -m(K_X+\Delta))$. Define 
$$S_m(E):=\max_{B_m}\ord_E(B_m),$$ 
where $B_m$ runs through all $m$-basis type divisors. It turns out 
$$\lim_m S_m(E)=S_{X, \Delta}(E), $$ 
e.g. \cite{BJ20}.
For the given log Fano pair $(X, \Delta)$, let $D$ be an effective $\bQ$-divisor on $X$ with $D\sim_\bQ -K_X-\Delta$. A simple computation implies the following relation
$$S_{X, \Delta+cD}(E)=(1-c)S_{X, \Delta}(E), $$
where $0\leq c< 1$ is a rational number and $E$ is any prime divisor over $X$.

\subsubsection{Pseudo-effective threshold}\label{subsubsec: peff}

Let $(X, \Delta)$ be a log Fano pair and $E$ a prime divisor over $X$. The \textit{pseudo-effective threshold} of $E$ with respect to $-K_X-\Delta$ is defined as follows:
$$T_{X,\Delta}(E):=\sup_{0\leq D\sim_\bQ -K_X-\Delta}\ord_E(D).$$
By \cite{BJ20}, we have the following relation between $S_{X, \Delta}(E)$ and $T_{X, \Delta}(E)$: 
$$\frac{d+1}{d}\leq \frac{T_{X,\Delta}(E)}{S_{X,\Delta}(E)}\leq d+1,$$ 
where $d$ is the dimension of $X$.

\subsection{K-stability}

Let $(X, \Delta)$ be a log Fano pair. Put 
$$\beta_{X,\Delta}(E):=A_{X,\Delta}(E)-S_{X,\Delta}(E).$$ 
By the works \cite{Fuj19, Li17}, one can define K-stability of a log Fano pair by beta criterion as follows.
\begin{definition}\label{def: kss}
Let $(X,\Delta)$ be a log Fano pair. 
We say that $(X,\Delta)$ is \emph{K-semistable} if $\beta_{X,\Delta}(E)\geq 0$ for any prime divisor $E$ over $X$.
\end{definition}

By the works \cite{Oda13, BHJ17}, one can define K-stability of a log Calabi-Yau pair by posing a singularity condition.

\begin{definition}
Let $(X, \Delta)$ be a log Calabi-Yau pair, i.e. $K_X+\Delta\sim_\bQ 0$. We say  $(X,\Delta)$ is \emph{K-semistable} if  $(X, \Delta)$ is log canonical.
\end{definition}

We also need the concept of delta invariant of a log Fano pair.
\begin{definition}
Let $(X, \Delta)$ be a log Fano pair, the \textit{delta invariant} of $(X, \Delta)$ is defined as 
$$\delta(X, \Delta):=\inf_{E} \frac{A_{X, \Delta}(E)}{S_{X, \Delta}(E)}, $$
where $E$ runs through all prime divisors over $X$ (e.g. \cite{FO18, BJ20}). By Definition \ref{def: kss} we see that $(X, \Delta)$ is K-semistable if and only if $\delta(X, \Delta)\geq 1$.
\end{definition}

We will encounter log pairs with $\bR$-coefficients in a special setting. 

\begin{definition}\label{def: real kss}
Let $X$ be a Fano variety and $0\leq D\sim_\bR -K_X$, we say $(X, cD)$ for some $0\leq c<1$ is a log Fano pair (with $\bR$-coefficients) if $(X, cD)$ is klt. In this setting, we also define the delta invariant of $(X, cD)$ as follows:
$$\delta(X, cD):=\inf_{E} \frac{A_{X, cD}(E)}{(1-c)S_{X}(E)}, $$
where $E$ runs through all prime divisors over $X$. We say that $(X, cD)$ is K-semistable if $\delta(X, cD)\geq 1$. We say that $(X, D)$ is K-semistable if $(X, D)$ is log canonical. 
\end{definition}
Notation as above, define 
$\tilde{\delta}(X, cD):=\min\{\delta(X, cD), 1\},$
then we have the following generalization of constructibility property for delta invariants in \cite{BLX22,Xu20}.

\begin{lemma}{\rm {(\cite{BLX22, Xu20})}}\label{lem: real constructibility}
Let $\pi: \mX\to B$ be a flat family of Fano varieties of dim $d$ such that $-K_{\mX/B}$ is a relatively ample $\bQ$-line bundle on $\mX$, and $B$ is a normal base. Suppose $\mD$ is an effective $\bR$-divisor on $\mX$ such that every component of $\mD$ is flat over $B$ and $\mD\sim_\bR -K_{\mX/B}$ over $B$. For any given real number $0<c<1$ such that $(\mX, c\mD)\to B$ is a family of log Fano pairs (with $\bR$-coefficients),  the set $\{\tilde{\delta}(\mX_t,c\mD_t)\ | \ t\in B\}$ is finite.
\end{lemma}

\begin{proof}
The idea is the same as \cite[Proposition 4.3]{BLX22}, though we are dealing with $\bR$-coefficients. Refer to \cite[Theorem 3.2]{Zhou23}.
\end{proof}

\subsection{Test configuration}

\begin{definition}\label{def: tc}
Let $(X,\Delta)$ be a log pair and $L$ an ample $\bQ$-line bundle on $X$. A \textit{test configuration} $\pi: (\mX,\Delta_\tc;\mL)\to \bA^1$ is a degenerating family over $\bA^1$ consisting of the following data:
\begin{enumerate}
\item $\pi: \mX\to \bA^1$ is a projective flat morphism from a normal variety $\mX$, $\Delta_\tc$ is an effective $\bQ$-divisor on $\mX$, and $\mL$ is a relatively ample $\bQ$-line bundle on $\mX$,
\item the family $\pi$ admits a $\bC^*$-action which lifts the natural $\bC^*$-action on $\bA^1$ such that $(\mX,\Delta_\tc; \mL)\times_{\bA^1}\bC^*$ is $\bC^*$-equivariantly isomorphic to $(X, \Delta; L)\times_{\bA^1}\bC^*$.
\end{enumerate}
\end{definition}

Suppose $(X,\Delta)$ is a log Fano pair and $L=-K_X-\Delta$. Let  $(\mX,\Delta_\tc; \mL)$ be a test configuration such that $\mL=-K_{\mX/\bA^1}-\Delta_\tc$. We call it a \textit{special test configuration} if  $(\mX, \mX_0+\Delta_{\tc})$ admits plt singularities (or equivalently, the central fiber $(\mX_0, \Delta_{\tc,0})$ is a log Fano pair). 
We will frequently use the following two results.

\begin{theorem}{\rm{(\cite{LXZ22})}}\label{thm: lxz22}
Let $(X, \Delta)$ be a log Fano pair with $\delta(X, \Delta)\leq 1$, then there exists a prime divisor $E$ over $X$ such that $\delta(X, \Delta)=\frac{A_{X, \Delta}(E)}{S_{X, \Delta}(E)}$.
\end{theorem}

\begin{theorem}{\rm{(\cite{BLZ22})}}\label{thm: blz22}
Let $(X, \Delta)$ be a log Fano pair with $\delta(X, \Delta)\leq 1$, and $E$ a prime divisor over $X$ such that $\delta(X, \Delta)=\frac{A_{X, \Delta}(E)}{S_{X, \Delta}(E)}$. Then $E$ induces a special test configuration $(\mX, \Delta_\tc)\to \bA^1$ of $(X, \Delta)$ such that $\delta(\mX_0, \Delta_{\tc, 0})=\delta(X, \Delta)$. In particular, if $(X, \Delta)$ is K-semistable, then the central fiber $(\mX_0, \Delta_{\tc, 0})$ is also K-semistable.
\end{theorem}

For a log pair $(X, \Delta)$ with $\bR$-coefficients, one could similarly define the test configuration for $(X, \Delta; L)$, where $L$ is a $\bQ$-ample line bundle on $X$. For example, let $X$ be a Fano variety and $D$ is an effective $\bR$-divisor on $X$ with $D\sim_\bR -K_X$,  we could define the test configuration of $(X, cD; -K_X)$ by the same way as in Definition \ref{def: tc}, where $0\leq c<1$ is a real number. Let $(\mX, c\mD; \mL)\to \bA^1$ be a test configuration of $(X, cD; -K_X)$, we say it is a special test configuration if $\mL=-K_{\mX/\bA^1}$ and $(\mX_0, c\mD_0)$ is a log Fano pair (with $\bR$-coefficients).

\begin{lemma}\label{lem: vanishing futaki}
 Let $X$ be a Fano variety and $D\sim_\bR -K_X$ an effective $\bR$-divisor on $X$. Suppose $0\leq c<1$ is a given real number such that $(X, cD)$ is a K-semistable log Fano pair (with $\bR$-coefficients), and $E$ is a prime divisor over $X$ such that it induces a special test configuration $(\mX, c\mD)$ of $(X, cD)$. If $(\mX_0, c\mD_0)$ is K-semistable, then 
 $$A_{X, cD}(E)-(1-c)S_{X}(E) =0.$$
 \end{lemma}

\begin{proof}
The proof is the same as that of \cite[Lemma 3.4]{Zhou23}.
\end{proof}

\subsection{Log bounded family}

Let $\mP$ be a set of projective varieties of dimension $d$. We say $\mP$ is bounded if there exists a projective morphism $\phi: \mY\to T$ between schemes of finite type such that for any $X\in \mP$, there exists a closed point $t\in T$ such that $X\cong \mY_t$. Let $\mP'$ be a set of log pairs with $\bR$-coefficients of dimension $d$, we say 
$\mP'$ is log bounded if there exist a projective morphism $\phi: \mY\to T$ between schemes of finite type and a reduced divisor $\mD$ on $\mY$ such that for any $(X, \Delta)\in \mP$, there exists a closed point $t\in T$ such that $(X, \red(\Delta))\cong (\mY_t, \mD_t)$. Here $\red(\Delta)$ means taking all the coefficients of components in $\Delta$ to be one. We have the following boundedness result for Fano varieties with K-stability.

\begin{theorem}{\rm{(\cite{Jiang20})}}\label{thm: jiang}
Fix a positive integer $d$ and two positive real numbers $v, \epsilon_0$. Then the following set lies in a bounded family:
$$\{X\ |\ \textit{$X$ is a Fano variety of dimension $d$ with $(-K_X)^d\geq v$ and $\delta(X)\geq \epsilon_0$}\}. $$
\end{theorem}

We also need the following lemma on log bounded family.

\begin{lemma}\label{lem: log bdd}
Let $f: (\mY, \sum_{j=1}^k\mD_j)\to T$ be a morphism such that 
\begin{enumerate}
\item $\mD_j$'s are effective $\bQ$-divisors on $\mY$;
\item $\mY$ is flat over $T$ and each component of $\mD_j$ is flat over $T$; 
\item for each closed point $t\in T$, the fiber $(\mY_t, \sum_{j=1}^k\mD_{k,t})$ is a log pair. 
\end{enumerate}
Let $\{(a_{i1},...,a_{ik})\}_i$ be a sequence of vectors tending to $(a_1,...,a_k)$. Suppose that $(\mY, \sum_{j=1}^ka_{ij}\mD_j)\to T$
is a family of log canonical pairs for every $i$ and $(\mY, \sum_{j=1}^ka_{j}\mD_j)\to T$ is a family of log pairs, then $(\mY, \sum_{j=1}^ka_{j}\mD_j)\to T$ is also a family of log canonical pairs.
\end{lemma}

\begin{proof}
 After a stratification of the base, we may assume $T$ is smooth and $f$ admits a fiberwise log resolution up to an \'etale base change. In this setting, the lemma just means: for a fixed log pair $(X, \sum_{j=1}^kD_j)$, if $(X, \sum_{j=1}^k a_{ij}D_j)$ is log canonical for every $i$ and $(X, \sum_{j=1}^k a_{j}D_j)$ is a log pair with $\bR$-coefficients, then $(X, \sum_{j=1}^k a_{j}D_j)$ is also log canonical. This is obvious.
\end{proof}

\subsection{Fano type variety}

Let $f: Y\to Z$ be a projective morphism such that $Y$ is a normal variety. We say that $Y$ is Fano type over $Z$ if there exists an effective $\bQ$-divisor $B$ on $Y$ such that $(Y, B)$ has klt singularities and $-K_Y-B$ is ample over $Z$. When $Z$ is a point, we just say that $Y$ is of Fano type or $Y$ is a Fano type variety.

A Fano type variety is naturally a Mori dream space (e.g. \cite{BCHM10}). This means that for a given projective morphism $Y\to Z$ such that $Y$ is Fano type over $Z$, one can run $D$-MMP over $Z$ for 
any $\bQ$-Cartier divisor $D$ to get a minimal model or a Mori fiber space.

\subsection{Complement}

\begin{definition}
Let $(X, \Delta)$ be a log pair. We say a $\bQ$-divisor $D\geq 0$ is a \emph{complement} of $(X, \Delta)$ if $(X, \Delta+D)$ is log canonical and $K_X+\Delta+D\sim_\bQ 0$; we say $D$ is an \emph{$N$-complement} for some positive integer $N$ if $D$ is a complement and $N(K_X+\Delta+D)\sim 0$.
\end{definition}

\begin{theorem}{\rm{(\cite{Birkar19})}}\label{thm: complement}
Let $(X, \Delta)$ be a log canonical pair with $-K_X-\Delta$ being nef and $X$ is a Fano type variety. Then there exists a positive number $N$ depending only on the dimension of $X$ and the coefficients of $\Delta$ such that $(X, \Delta)$ admits an $N$-complement.
\end{theorem}

\section{Convex sets and polytopes}\label{sec: polytope}

In this section, we introduce some notation and concepts on convex sets and polytopes, which we will use frequently.

\subsection{Convex sets and Polytopes}
\begin{definition}
Let $C$ be a closed subset of $\bR^k$,  we say that $C$ is a convex set (or $C$ is convex) if for any two points $\vec{a}, \vec{b}\in C$, the segment connecting $\vec{a}$ and $\vec{b}$ is also contained in $C$. 
View $C$ as a subset of $1\times \bR^k\subset \bR\times \bR^k$. Let $\vec{0}\in \bR^{k+1}$ be the origin. Denote by $\Cone(\vec{0}, C)$ the cone over $C$ which contains every line starting from $\vec{0}$ and passing through some point in $C$. We say that $\vec{a}\in C$ is an extreme point of $C$ if there is an extremal ray of $\Cone(\vec{0}, C)$ passing through $\vec{a}$.
\end{definition}

Let $H$ be a hyperplane in $\bR^k$ determined by the linear form $\sum_{j=1}^ka_jx_j-a=0$. 
We use the following notation:
$$H^+:=\{(x_1,...,x_k)\in \bR^k\ |\ \textit{$\sum_{j=1}^ka_jx_j-a>0$}\}, $$
$$H^{\geq 0}:=\{(x_1,...,x_k)\in \bR^k\ |\ \textit{$\sum_{j=1}^ka_jx_j-a\geq 0$}\}, $$
$$H^-:=\{(x_1,...,x_k)\in \bR^k\ |\ \textit{$\sum_{j=1}^ka_jx_j-a<0$}\}, $$
$$H^{\leq 0}:=\{(x_1,...,x_k)\in \bR^k\ |\ \textit{$\sum_{j=1}^ka_jx_j-a\leq 0$}\}. $$

\begin{definition}\label{def: separation}
Let $C\subset \bR^k$ be a convex set and $H\subset \bR^k$ a hyperplane determined by the form $\sum_{j=1}^ka_jx_j-a=0$. We say that $H$ separates $C$ if both $H^+\cap C$ and $H^-\cap C$ are non-empty sets.
\end{definition}

\begin{definition}\label{def: polytope}
Let $P\subset \bR^k$ be a compact subset. We say that $P$ is a polytope if $P$ is of the following form of finite intersection:
$$P=\bigcap_{i=1}^s H_i^{\geq 0}, $$
where $H_i$ ($i=1,...,s$) are hyperplanes in $\bR^k$. We say $P$ is a rational polytope if every extreme point of $P$ is rational in $\bR^k$ (i.e. all the coordinates of the point are rational numbers). We say the subset $F\subset P$ is a face of $P$ if $F$ is of the following form of finite intersection
$$F=\bigcap_{j=1}^rH_{i_{j}}\cap P $$
for some $\{i_1,i_2,...,i_r\}\subset\{1,2,...,s\}$. Suppose $\dim P=l$, we say that a face of dimension $l-1$ is a facet of $P$. 
\end{definition}

\begin{definition}
Let $C\subset \bR^k$ be a convex set and $\vec{w}$ an extreme point of $C$. We say that $C$ is locally a polytope at $\vec{w}$ if there exists an open neighborhood of $\vec{w}$, denoted by $U_{\vec{w}}$, such that 
$$C\cap U_{\vec{w}}=\bigcap_{i=1}^s H_i^{\geq 0}\cap U_{\vec{w}}$$ 
for some hyperplanes $H_i, \ i=1,...,s$. 
\end{definition}

\begin{definition}
Let $Z\subset \bR^k$ be an affine subspace, i.e. $Z=\cap_{i=1}^s H_i$ for some hyperplanes $H_i$. We say $Z$ is rational if the rational points on $Z$ are dense.
\end{definition}

We have the following characterization of rationality of affine subspaces.

\begin{lemma}\label{lem: rational hyperplane}
Let $H\subset \bR^k$ be a hyperplane, then $H$ is rational if and only if $H$ can be formulated by $\sum_{j=1}^ka_jx_j-a=0$, where $a_j$ and  $a$ are all rational numbers.
\end{lemma}

\begin{proof}
Suppose $H$ could be formulated by $\sum_{j=1}^ka_jx_j-a=0$, where $a_j$ and  $a$ are all rational numbers,  then clearly $H$ is rational. Conversely, suppose $H$ is rational. If $a\neq 0$, we may assume $a=1$, and there exist $k$ rational points on $H$, denoted by $\vec{w}_i:=(x_{i1},...,x_{ik}),\ i=1,...,k$, such that the rational matrix $(x_{ij})_{1\leq i,j\leq k}$ is invertible. Hence the $k$ equations 
$$\sum_{j=1}^ka_jx_{ij}=1,\ i=1,...,k$$ 
give a unique rational solution $(a_1,...,a_k)$. If $a=0$, then there exist $k-1$ rational points on $H$, denoted by $\vec{w}_i:=(x_{i1},...,x_{ik}),\ i=1,...,k-1$, such that the rational matrix $(x_{ij})_{1\leq i\leq k-1, 1\leq j\leq k}$ is of rank $k-1$. We may assume $a_k=1$ and the matrix $(x_{ij})_{1\leq i,j\leq k-1}$ is invertible, then the $k-1$ equations 
$$\sum_{j=1}^{k-1}x_{ij}a_j+x_{ik}=0,\ i=1,...,k-1$$ 
admit a unique rational solution $(b_1,...,b_{k-1})$. We are done by replacing $(a_1,...,a_k)$ with $(b_1,...,b_{k-1}, 1)$.
\end{proof}

\begin{lemma}\label{lem: rational subspace}
Let $Z\subset \bR^k$ be an affine subspace, then $Z$ is rational if and only if Z could be formulated as the intersection of rational hyperplanes.
\end{lemma}

\begin{proof}
Suppose $\dim\ Z=k-s$ and $Z$ could be formulated by the intersection
$Z=\cap_{i=1}^s H_i, $
where $H_i$ is determined by $\sum_{j=1}^ka_{ij}x_j-a_i=0,\ i=1,...,s$.

Assume $Z$ passes through the origin, i.e. $a_i=0$ for every $i$. We first suppose all $a_{ij}$ ($1\leq i\leq s, 1\leq j\leq k$) are rational numbers, then up to a linear transformation of $\bR^k$ via a rational invertible matrix, we easily see that $Z$ is rational. Conversely, suppose $Z$ is rational, then there exist $k-s$ rational points on $Z$, denoted by $\vec{w}_i:=(x_{i1},...,x_{ik}),\ i=1,...,k-s$, such that the matrix $(x_{ij})_{1\leq i\leq k-s, 1\leq j\leq k}$ is of rank $k-s$. Thus the $s$-dimensional linear subspace determined by 
$$\sum_{j=1}^k y_j x_{ij}=0,\ i=1,...,k-s, $$
is rational, and there exist $s$ rational points $\vec{v}_i:=(y_{i1},...,y_{ik}),\ i=1,...,s$, such that they span the $s$-dimensional linear subspace. Replacing the matrix $(a_{ij})_{1\leq i\leq s, 1\leq j\leq k}$ with $(y_{ij})_{1\leq i\leq s, 1\leq j\leq k}$, we see that $Z$ is the intersection of rational hyperplanes.

Assume $Z$ does not pass through the origin. Suppose all $a_{ij}, a_i$ ($1\leq i\leq s, 1\leq j\leq k$) are rational numbers, then one could easily find a rational point $(b_1,...,b_k)\in Z$ such that $\sum_{j=1}^ka_{ij}(x_j-b_j)=0$ for every $i$. Via transforming $Z$ along the vector $(b_1,...,b_k)$, we get a rational linear subspace, which conversely implies that $Z$ is rational.
Conversely, we suppose $Z$ is rational, then there exists a rational point $(b_1,...,b_k)$ on $Z$ such that $\sum_{j=1}^ka_{ij}(x_j-b_j)=0$ for every $i$. Via transforming $Z$ along the vector $(b_1,...,b_k)$, we get a rational linear subspace.
By what we have proved, one could replace $a_{ij}$ with rational numbers while preserving the space $Z$, and $a_i=\sum_{j=1}^ka_{ij}b_j$ is  then also rational.
\end{proof}

\begin{corollary}\label{cor: rational vertex}
Let $P\subset \bR^k$ be a polytope. Suppose $\vec{w}$ is an extreme point and $\vec{w}$ is given by the intersection of facets $F_1,..., F_s$, i.e. $\cap_{i=1}^s F_s=\vec{w}$. If the rational points on $F_i$ for $\ i=1,...,s$ are dense, then $\vec{w}$ is rational.
\end{corollary}

\begin{proof}
We denote by $Z_i$ the minimal affine subspace containing $F_i$, then clearly $Z_i$ is rational for very $i$. By Lemma \ref{lem: rational subspace}, $Z_i$ could be formulated as the intersection of rational hyperplanes. Thus $\vec{w}$ could be formulated as the intersection of some rational hyperplanes $H_l, l=1,...r$, where $H_l$ is determined by $\sum_{j=1}^ka_{lj}x_j-a_l=0$, and these equations admit the unique solution $\vec{w}$. By Lemma \ref{lem: rational hyperplane}, we could assume all $a_{lj}, a_l$ are rational numbers, thus the solution $\vec{w}$ is also rational.
\end{proof}

We also need the following local characterization of a polytope.

\begin{lemma}\label{lem: local char polytope}
Let $C\subset \bR^k$ be a compact convex set and $\vec{w}\in C$ an extreme point of $C$. Suppose $C$ is not locally a polytope at $\vec{w}$, then there exists a sequence of points $\{\vec{w}_i\}_{i=1}^\infty$ such that
\begin{enumerate}
\item $\vec{w}_i$ is on the boundary of $C$ for every $i$;
\item $\vec{w}_i$ tends to a point on the boundary of $C$, denoted by $\vec{w}'$ (not necessarily $\vec{w}$);
\item any hyperplane passing through $\vec{w}_i$ and $\vec{w}'$ separates $C$ for every $i$.
\end{enumerate}
\end{lemma}

\begin{proof}
We denote $\partial C$ to be the boundary of $C$.
Suppose for any $0<\epsilon\ll 1$, there exists a point $\vec{v}\in \partial C$ with $||\vec{w}-\vec{v}||<\epsilon$ such that the segment connecting $\vec{w}$ and $\vec{v}$ does not lie on  $\partial C$, then one could easily find a sequence $\vec{w}_i$ on $\partial C$ such that $\vec{w}_i$ tends to $\vec{w}$ and any hyperplane passing through $\vec{w_i}$ and $\vec{w}$ separates $C$.

Now we suppose there exists a positive number $0<\epsilon<1$ such that for any point $\vec{v}\in \partial C$ with $||\vec{w}-\vec{v}||<\epsilon$, the segment connecting $\vec{w}$ and $\vec{v}$ lies on $\partial C$. In this case, there exist finitely many points on $\partial C$, denoted by $\vec{v}_l,\ l=1,...,s,$ such that all these $\vec{v}_l$ lie on the same hyperplane $H$, and $H\cap C$ is not a polytope (otherwise $C$ is locally a polytope at $\vec{w}$).  For the convex set $H\cap C$, 
suppose it is locally a polytope at every extreme point, then we can find a sequence of extreme points $\vec{w}_i$ of $H\cap C$ such that $\vec{w}_i$ tends to a point on $\partial (H\cap C)$, denoted by $\vec{w}'$, and any hyperplane passing through $\vec{w_i}$ and $\vec{w}'$ separates $H\cap C$, hence also separates $C$.  Thus we may assume that there exists  an extreme  point $\vec{w}^{(1)}$ of $H\cap C$ such that $H\cap C$ is not locally a polytope at $\vec{w}^{(1)}$. Replace $(C, \vec{w})$ with $(H\cap C, \vec{w}^{(1)})$, one could apply the same analysis as before.  It then suffices to show the argument for the case $\dim C=2$. In this case, it is clear that for any $0<\epsilon\ll 1$, there exists a point $\vec{v}\in \partial C$ with $||\vec{w}-\vec{v}||<\epsilon$ such that the segment connecting $\vec{w}$ and $\vec{v}$ does not lie on  $\partial C$, thus we get the desired sequence.
\end{proof}

\begin{remark}
In Lemma \ref{lem: local char polytope}, if the rational points are dense on $\partial C$, then one can choose a rational sequence $\{\vec{w}_i\}_i$.
\end{remark}

\subsection{Log canonical polytope}

Let $(X, \sum_{j=1}^k\Delta_j)$ be a log pair, where $X$ is log canonical and $\Delta_j, \ j=1,...,k$, are effective $\bQ$-Cartier divisors on $X$. Then we define the following set
\begin{align*}
 \LC(X, \sum_{j=1}^k\Delta_j)
:=\ \{(a_1,...,a_k)\in \bR_{\geq 0}^k\ |\ \textit{$(X, \sum_{j=1}^ka_j\Delta_j) $ is log canonical\ }\}. 
\end{align*}

\begin{lemma}\label{lem: lc polytope}
Let $(X, \sum_{j=1}^k\Delta_j)$ be a log pair, where $X$ is log canonical and $\Delta_j, \ j=1,...,k$, are effective $\bQ$-Cartier divisors on $X$, then $\LC(X, \sum_{j=1}^k\Delta_j)$ is a rational polytope (see Definition \ref{def: polytope}).
\end{lemma}

\begin{proof}
Let $f: Y\to (X, \sum_{j=1}^k\Delta_j)$ be a log resolution and write
$$K_Y=f^*(K_X+\sum_{j=1}^kx_j\Delta_j)+\sum_{i} l_i(x_1,...,x_k) E_i, $$
where $E_i$ are smooth prime divisors on $Y$ and $l_i(x_1,...,x_k)$ are linear functions on $x_1,...,x_k$ with $\bQ$-coefficients. Thus we see that $\LC(X, \sum_{j=1}^k\Delta_j)$ is given by the equations 
$$l_i(x_1,...,x_k)\geq -1,$$ 
which produces a rational polytope by Lemma \ref{lem: rational subspace}.
\end{proof}

\section{Boundedness}\label{sec: boundedness}

In \cite{LZ23}, we proved that the set $\mE(d,k,v,I)$ is log bounded. This will be a crucial ingredient for us to finally confirm the finite chamber decomposition of $\overline{\Delta^{k}}$ as stated in Theorem \ref{thm: main1} and Corollary \ref{cor: chamber decomposition}. For the readers' convenience, we include a complete proof of the boundedness here.

\begin{theorem}\label{thm: boundedness}
Fix two positive integers $d$ and $k$, a positive number $v$,  and a positive integer $I$, then the set $\mE:=\mE(d,k,v, I)$ is log bounded.
\end{theorem}

\begin{proof}
We divide the proof into several steps.

\

\textit{Step 1}. In this step, we explain that it is enough to show that the set 
$$\mG:=\{X\ |\ \textit{$(X, \sum_{j=1}^kD_j)\in \mE$}\}$$
is bounded. Suppose $\mG$ is bounded, then for each $(X, \sum_{j=1}^kD_j)\in \mE$, there exists a very ample line bundle $L_X$ on $X$ such that 
$$L_X^d\leq M(d)\quad \text{and} \quad -K_X\cdot L_X^{d-1}\leq N(d),$$ 
where $M(d)$ and $N(d)$ are positive numbers depending only on the dimension $d$. For each $j$, it is clear that $D_j\cdot L_X^{d-1}\leq N(d)$. As the coefficients of $D_j$ are contained in $\frac{1}{I}\bN$, we see that the degree of each component of $D_j$ is upper bounded. By applying a standard argument on Chow scheme we see that $D_j$ lies in a bounded family, hence $\mE$ is log bounded. The rest of the proof is devoted to the boundedness of $\mG$.

\

\textit{Step 2}. In this step, we want to replace $(X, \sum_{j=1}^kD_j)$ with $(X, \sum_{j=1}^kD_j')$ such that  $(X, \sum_{j=1}^kD_j')\in \mE$ and $(X, D_j')$ is a  log canonical Calabi-Yau pair for every $1\leq j\leq k$. By Theorem \ref{thm: complement}, there exists a positive integer $m(d)$ depending only on the dimension $d$ such that $X$ admits $m(d)$-complements. We may choose $m(d)$ properly such that $I|m(d)$. Put 
$$\bP:=\frac{1}{m(d)}|-m(d)K_X|.$$ 
Then we see that $D_j\in \bP$ for every $1\leq j\leq k$. Since $(X, \sum_{j=1}^kD_j)\in \mE$, there exist some rational numbers $0\leq c_j<1$ such that $(X, \sum_{j=1}^k c_jD_j)$ is K-semistable. We are ready to replace $D_1$ with some $D_1'$.
Let $\mD\subset X\times \bP$ be the universal divisor associated to the linear system $\frac{1}{m(d)}|-m(d)K_X|$ and consider the universal family
$$ (X\times \bP, c_1\mD+\sum_{j=2}^kc_jD_j\times \bP)\to \bP .$$
Since there exists a fiber which is K-semistable, i.e. $(X, \sum_{j=1}^kc_jD_j$), one could find an open subset $U\subset \bP$ such that $(X, c_1\mD_t+\sum_{j=2}^kc_jD_j)$ is K-semistable for any $t\in U$ (see \cite[Proposition 4.3]{BLX22}). Recall that $X$ admits $m(d)$-complements, so we conclude that $(X, \mD_t)$ is log canonical for general $t\in U$. Replacing $D_1$ with such a $\mD_t$ and replacing $I$ with a bounded multiple, we obtain $D_1'$ as required. By the same way, we could replace other $D_j, 2\leq j\leq k$, step by step.

\

\textit{Step 3}. By step 2, we assume $(X, \sum_{j=1}^kD_j)\in \mE$ satisfies that $(X, D_j)$ is log canonical for every $1\leq j\leq k$. For such log pair $(X, \sum_{j=1}^kD_j)$, we define the following invariant:
$$\mu(X, \sum_{j=1}^kD_j):=\inf\Bigg\{\sum_{j=1}^kc_j\ |\ \textit{$(X, \sum_{j=1}^kc_jD_j)$ is K-semistable}\Bigg\}. $$
It is clear that $0\leq \mu(X, \sum_{j=1}^kD_j)<1$. We aim to show that there is a gap between $\mu(X, \sum_{j=1}^kD_j)$ and $1$. More precisely, there exists a positive number $0<\epsilon_0(d,k,v, I)<1$ depending only on $d,k,v,I$ such that
$$1-\mu(X, \sum_{j=1}^kD_j)\geq \epsilon_0. $$
Suppose not, we could find a sequence of log pairs $(X_i, \sum_{j=1}^kD_{ij})$ satisfying the following conditions:
\begin{enumerate}
\item $(X_i, \sum_{j=1}^kD_{ij})\in \mE$ for every $i$,
\item $(X_i, D_{ij})$ is log canonical for every $1\leq j\leq k$ and $i$,
\item $\mu_i:=\mu(X_i, \sum_{j=1}^kD_{ij})$ is an increasing sequence tending to $1$.
\end{enumerate}
By the definition of $\mu_i$, there exists an increasing sequence of rational numbers $a_i<\mu_i$ tending to $1$ such that
$$(X_i, \sum_{j=1}^k\frac{a_i}{k}D_{ij}) $$
is K-unstable for every $i$. By Theorem \ref{thm: lxz22} and Theorem \ref{thm: blz22}, for each $(X_i, \sum_{j=1}^k\frac{a_i}{k}D_{ij}) $, there exists a prime divisor $E_i$ over $X_i$ such that $E_i$ computes the delta invariant of $(X_i, \sum_{j=1}^k\frac{a_i}{k}D_{ij})$ and $E_i$ induces a special test configuration 
$$(\mX_i, \sum_{j=1}^k\frac{a_i}{k}\mD_{ij})\to \bA^1.$$
Subtracting those $a_i$ which are not close enough to $1$, we may assume that the central fiber of the test configuration (after changing the coefficients), i.e. $(\mX_{i, 0}, \sum_{j=1}^k\frac{1}{k}\mD_{ij, 0})$, is log canonical (e.g. \cite{HMX14}). This means that the test configuration degenerates the log canonical Calabi-Yau pair $(X_i, \sum_{j=1}^k\frac{1}{k}D_{ij})$ to another log canonical Calabi-Yau pair $(\mX_{i, 0}, \sum_{j=1}^k\frac{1}{k}\mD_{ij, 0})$.
By \cite[Lemma 2.8]{Zhou23}, $E_i$ is an lc place of $(X_i, \sum_{j=1}^k\frac{1}{k}D_{ij})$, which forces that $E_i$ is an lc place of $(X_i, D_{ij})$ for every $1\leq j\leq k$. Recall that $E_i$ computes the delta invariant of the K-unstable log Fano pair $(X_i, \sum_{j=1}^k\frac{a_i}{k}D_{ij})$. We have
\begin{align*}
\beta_{X_i, \sum_{j=1}^k\frac{a_i}{k}D_{ij}}(E_i)\ &=\ A_{X_i, \sum_{j=1}^k\frac{a_i}{k}D_{ij}}(E_i) -S_{X_i, \sum_{j=1}^k\frac{a_i}{k}D_{ij}}(E_i)\\
&=\ (1-a_i)\{A_{X_i}(E_i)-S_{X_i}(E_i)\}\\
&<\ 0.
\end{align*}
On the other hand, $(X_i, \sum_{j=1}^kc_{ij}D_{ij})$ is a K-semistable log Fano pair  for some rational $0\leq c_{ij}<1$, therefore,
\begin{align*}
\beta_{X_i, \sum_{j=1}^kc_{ij}D_{ij}}(E_i)\ &=\ A_{X_i, \sum_{j=1}^kc_{ij}D_{ij}}(E_i) -S_{X_i, \sum_{j=1}^kc_{ij}D_{ij}}(E_i)\\
&=\ (1-\sum_{j=1}^kc_{ij})\{A_{X_i}(E_i)-S_{X_i}(E_i)\}\\
&\geq\ 0,
\end{align*}
which is a contradiction. This contradiction implies the existence of the gap we want.

\

\textit{Step 4}. Combining step 2 and step 3, we see that for each $(X, \sum_{j=1}^kD_j)\in \mE$, one can find another log pair $(X, \sum_{j=1}^kD_j')\in \mE$ such that $(X, \sum_{j=1}^kc_jD_j')$ is K-semistable for some numbers $0\leq c_j<1$ with
$$1-\sum_{j=1}^kc_j\geq \epsilon_0.$$
Thus we have
\begin{align*}
\frac{A_{X}(E)}{(1-\sum_{j=1}^kc_j)S_X(E)}\ &\geq  \ \frac{A_{X, \sum_{j=1}^kc_jD'_j}(E)}{S_{X, \sum_{j=1}^kc_jD'_j}(E)}\geq 1
\end{align*}
for any prime divisor $E$ over $X$. Hence,
$$\frac{A_X(E)}{S_X(E)}\geq 1-\sum_{j=1}^kc_j\geq \epsilon_0, $$
for any prime divisor $E$ over $X$. This says that the delta invariant of $X$ is bounded from below by $\epsilon_0(d, I)$. By Theorem \ref{thm: jiang}, the set $\mG$ defined in step 1 lies in a bounded family. The proof is finished.
\end{proof}

\section{On the shape of the K-semistable domain}\label{sec: rational polytope}

Fix two positive integers $d$ and $k$, a positive number $v$, and a positive integer $I$. In this section, we will prove several properties of K-semistable domains of log pairs in $\mE:=\mE(d,k,v,I)$. In particular, we will show that they are rational polytopes. 

\begin{lemma}\label{lem: dense lemma}
Let $X$ be a Fano variety of dimension $d$, and $B_1, B_2$ are two effective $\bQ$-divisors on $X$ such that 
$B_i$'s are proportional to $-K_X$. Suppose $(X, tB_1+(1-t)B_2)$ is a K-semistable log Fano pair for any $t\in (t_0, t_0+\epsilon)\cap \bQ$ {\rm{(resp. $t\in (t_0-\epsilon, t_0)\cap \bQ$)}} but a K-unstable log Fano pair for any $t\in (t_0-\epsilon, t_0)\cap\bQ$ {\rm{(resp. $t\in (t_0, t_0+\epsilon)\cap \bQ$)}}, where $0<t_0<1$ and $0<\epsilon\ll 1$.  Then $t_0$ is rational.
\end{lemma}

\begin{proof}

The idea is essentially the same as that of \cite[Lemma 4.3]{Zhou23}. We present the proof for the readers' convenience.
We assume $(X, tB_1+(1-t)B_2)$ is a K-semistable log Fano pair for any $t\in (t_0, t_0+\epsilon)\cap \bQ$ but a K-unstable log Fano pair for any $t\in (t_0-\epsilon, t_0)\cap\bQ$. The other case can be proven similarly. 

Choosing a strictly increasing sequence of positive rational numbers $\{a_i\}_{i=1}^\infty$ tending to $t_0$, then $(X, a_iB_1+(1-a_i)B_2)$ is not K-semistable for any $i$. By Theorem \ref{thm: lxz22} and Theorem \ref{thm: blz22}, there exists a prime divisor $E_i$ over $X$ such that 
$$A_{X, a_iB_1+(1-a_i)B_2}(E_i)-\delta(X, a_iB_1+(1-a_i)B_2)\cdot S_{X, a_iB_1+(1-a_i)B_2}(E_i)=0, $$
and $E_i$ induces a special test configuration of $(X, a_iB_1+(1-a_i)B_2)$, denoted by
$$(\mX_i, a_i\mB_{i1}+(1-a_i)\mB_{i2})\to \bA^1,$$
such that 
$$\delta(\mX_{i,0}, a_i\mB_{i1,0}+(1-a_i)\mB_{i2,0})=\delta(X, a_iB_1+(1-a_i)B_2)<1. $$
We first show that the set $\{(\mX_{i,0}, \mB_{i1,0}+\mB_{i2,0})\}$ lies in a log bounded family. Suppose 
$$B_1\sim_\bQ -b_1K_X\quad \text{and}\quad B_2\sim_\bQ -b_2K_X.$$
It is clear  that $\vol(-K_{\mX_{i,0}})=\vol(-K_X)$ and 
$$\frac{A_{\mX_{i,0}}(E)}{S_{\mX_{i,0}}(E)}\geq  \{1-b_1a_i-b_2(1-a_i)\}\cdot\frac{A_{\mX_{i,0}, a_i\mB_{i1,0}+(1-a_i)\mB_{i2,0}}(E)}{S_{\mX_{i,0}, a_i\mB_{i1,0}+(1-a_i)\mB_{i2,0}}(E)}$$
for every prime divisor $E$ over $\mX_{i,0}$. Denote by 
$$c_i:=b_1a_i+b_2(1-a_i),$$ 
which tends to 
$$c:=b_1t_0+b_2(1-t_0).$$ 
Then we have the following estimate $(\clubsuit)$
\begin{align*}
\delta(\mX_{i,0}) &\ \geq (1-c_i)\cdot\delta(\mX_{i,0}, a_i\mB_{i1,0}+(1-a_i)\mB_{i2,0})\\
&\ =(1-c_i)\cdot\delta(X,a_iB_1+(1-a_i)B_2).
\end{align*}
Note that $\delta(X,t_0B_1+(1-t_0)B_2)\geq 1$ in the sense of Definition \ref{def: real kss}, thus 
$$\frac{A_{X, t_0B_1+(1-t_0)B_2}(F)}{(1-c)S_X(F)}\geq 1$$ for any prime divisor $F$ over $X$. We have the following estimate $(\spadesuit)$ 
\begin{align*}
&\ \delta(X,a_iB_1+(1-a_i)B_2)\\
=&\ \inf_F \frac{A_{X,a_iB_1+(1-a_i)B_2}(F)}{S_{X,a_iB_1+(1-a_i)B_2}(F)}\\
=&\ \inf_F \frac{A_{X,t_0B_1+(1-t_0)B_2}(F)+(t_0-a_i)\ord_F(B_1)+(a_i-t_0)\ord_F(B_2)}{\frac{1-c_i}{1-c}(1-c)S_X(F)}\\
\geq&\ \inf_F \frac{A_{X,t_0B_1+(1-t_0)B_2}(F)-(t_0-a_i)(b_1+b_2)T_X(F)}{\frac{1-c_i}{1-c}(1-c)S_X(F)}\\
\geq&\   \frac{1-c}{1-c_i}- \frac{(d+1)(t_0-a_i)(b_1+b_2)}{1-c_i},
\end{align*}
where $F$ runs through all prime divisors over $X$. Note that the first inequality follows from the definition of the pseudo-effective threshold and the second inequality follows from the estimate of $\frac{T_X(F)}{S_X(F)}$ (see Section \ref{subsubsec: peff}).
Combining $(\clubsuit)$ and $(\spadesuit)$, we see that $\delta(\mX_{i,0})$ admits a positive lower bound which does not depend on $i$ since 
$$\delta(\mX_{i,0})\geq 1-c-(d+1)(t_0-a_i)(b_1+b_2), $$
and $a_i$ tends to $t_0$. 
Thus the set $\{\mX_{i,0}\}$ lies in a bounded family by Theorem \ref{thm: jiang}, which implies that the set $\{(\mX_{i,0},\mB_{i1,0}+\mB_{i2,0})\}$ lies in a log bounded family as 
$$\mB_{ij,0}\sim_\bQ-b_jK_{\mX_{i,0}},\ \ j=1,2, $$
and the Weil indices of $\mB_{ij,0}$ ($j=1,2$) are upper bounded. We now consider the set $\{(\mX_{i,0}, t_0\mB_{i1,0}+(1-t_0)\mB_{i2,0})\}$, which lies in a log bounded family. We have the following estimate $(\bigstar)$ 
\begin{align*}
&\ \delta(\mX_{i,0},t_0\mB_{i1,0}+(1-t_0)\mB_{i2,0})\\
= &\ \inf_{F}\frac{A_{\mX_{i,0},t_0\mB_{i1,0}+(1-t_0)\mB_{i2,0}}(F)}{(1-c)S_{\mX_{i,0}}(F)}\\
= &\ \inf_F\frac{A_{\mX_{i,0},a_i\mB_{i1,0}+(1-a_i)\mB_{i2,0}}(F)+(a_i-t_0)\ord_F(\mB_{i1,0})+(t_0-a_i)\ord_F(\mB_{i2,0})}{\frac{1-c}{1-c_i}S_{\mX_{i,0},a_i\mB_{i1,0}+(1-a_i)\mB_{i2,0}}(F)}\\
\geq&\ \inf_F \frac{1-c_i}{1-c}\cdot\frac{A_{\mX_{i,0},a_i\mB_{i1,0}+(1-a_i)\mB_{i2,0}}(F)}{S_{\mX_{i,0},a_i\mB_{i1,0}+(1-a_i)\mB_{i2,0}}(F)}-\frac{(t_0-a_i)(d+1)(b_1+b_2)}{1-c}\\
= &\ \frac{1-c_i}{1-c}\cdot\delta(\mX_{i,0}, a_i\mB_{i1,0}+(1-a_i)\mB_{i2,0})-\frac{(t_0-a_i)(d+1)(b_1+b_2)}{1-c},
\end{align*}
where $F$ runs through all prime divisors over $\mX_{i,0}$.
Combining the estimate $(\spadesuit)$ of 
$$\delta(\mX_{i,0}, a_i\mB_{i1,0}+(1-a_i)\mB_{i2,0})=\delta(X, a_iB_1+(1-a_i)B_2)$$
with the estimate $(\bigstar)$, we see the following
$$\delta(\mX_{i,0}, t_0\mB_{i1,0}+(1-t_0)\mB_{i2,0})\geq 1-  \frac{2(d+1)(t_0-a_i)(b_1+b_2)}{1-c}.$$
Since $\{(\mX_{i,0}, t_0\mB_{i1,0}+(1-t_0)\mB_{i2,0})\}$ is log bounded, by applying Lemma \ref{lem: real constructibility}, we see that $\delta(\mX_{i,0}, t_0\mB_{i1,0}+(1-t_0)\mB_{i2,0})\geq 1$ for $i\gg1$. By Lemma \ref{lem: vanishing futaki}, this implies 
\begin{align*}
&\ A_{X,t_0B_1+(1-t_0)B_2}(E_i)-S_{X, t_0B_1+(1-t_0)B_2}(E_i)\\
=&\  t_0\{-\ord_{E_i}(B_1)+\ord_{E_i}(B_2)+b_1S_X(E_i)-b_2S_X(E_i)\}\\
&\ + \{A_X(E_i)-\ord_{E_i}(B_2)-S_X(E_i)+b_2S_X(E_i)\}\\
=&\ 0.
\end{align*}
By the proof of \cite[Theorem 4,7]{BLZ22} we know that $S_X(E_i)$ is a rational number. Combining the inequality
$$A_{X,a_iB_1+(1-a_i)B_2}(E_i)-S_{X, a_iB_1+(1-a_i)B_2}(E_i)<0,$$
we see that 
$$\ord_{E_i}(B_1)-\ord_{E_i}(B_2)-b_1S_X(E_i)+b_2S_X(E_i)\ne 0. $$
Thus 
$$t_0=\frac{A_X(E_i)-\ord_{E_i}(B_2)-S_X(E_i)+b_2S_X(E_i)}{\ord_{E_i}(B_1)-\ord_{E_i}(B_2)-b_1S_X(E_i)+b_2S_X(E_i)} $$
is a rational number.
\end{proof}

From Lemma \ref{lem: dense lemma} we directly see that there are dense rational points on the boundary of $\Kss(X, \sum_{j=1}^kD_j)$, where $(X, \sum_{j=1}^kD_j)\in \mE$.

\begin{lemma}\label{lem: rationality lemma}
Let $(X, \sum_{j=1}^kD_j)$ be a log pair in $\mE:=\mE(d,k,v,I)$ and $\vec{w}\in \overline{\Delta^k}$ is an extreme point of $\Kss(X, \sum_{j=1}^kD_j)$. If $\Kss(X, \sum_{j=1}^kD_j)$ 
is locally a polytope at $\vec{w}$, then $\vec{w}$ is rational.
\end{lemma}

\begin{proof}
Since the K-semistable domain $\Kss(X, \sum_{j=1}^kD_j)$ is locally a polytope at $\vec{w}$, we may assume that  $\Kss(X, \sum_{j=1}^kD_j)$ is cut out by hyperplanes $H_1, H_2,..., H_s$ locally at $\vec{w}$, and the intersection of these hyperplanes produces the vertex $\vec{w}$. The facet at $\vec{w}$ is then of the form $H_{i_1}\cap...\cap H_{i_l}$ for some $1\leq i_1<...<i_l\leq s$. We denote the facets at $\vec{w}$ by $F_1,...,F_m$, then $\vec{w}$ is the vertex of the intersection of these facets. 
To show that $\vec{w}$ is rational, it suffices to show that every facet $F_i$ is rational (see Corollary \ref{cor: rational vertex}). By Lemma \ref{lem: dense lemma}, we know that the rational points on $F_i\cap \Kss(X, \sum_{j=1}^kD_j)$ are dense, which implies that $F_i$ is rational. Since $\vec{w}$ is the vertex of the intersection of $F_i$ ($i=1,...,m$), $\vec{w}$ is also rational.
\end{proof}

\begin{lemma}\label{lem: polytope lemma}
Let $(X, \sum_{j=1}^kD_j)$ be a log pair in $\mE:=\mE(d,k,v,I)$ and $\vec{w}\in \overline{\Delta^k}$ is an extreme point of $\Kss(X, \sum_{j=1}^kD_j)$. Then $\Kss(X, \sum_{j=1}^kD_j)$ 
is locally a polytope at $\vec{w}$.
\end{lemma}

\begin{proof}
Note that $\Kss(X, \sum_{j=1}^kD_j)$ is a convex set. If $\Kss(X, \sum_{j=1}^kD_j)$ is contained in the hyperplane $\{\sum_{j=1}^k x_j=1\}$, then $\Kss(X, \sum_{j=1}^kD_j)$ is a rational polytope by Lemma \ref{lem: lc polytope}.  Thus we may assume $\Kss(X, \sum_{j=1}^kD_j)$ is not contained in $\{\sum_{j=1}^k x_j=1\}$.
Suppose $\Kss(X, \sum_{j=1}^kD_j)$ is not locally a polytope at $\vec{w}$, by Lemma \ref{lem: local char polytope}, there exists a sequence of rational points $\{\vec{w}_i\}_{i=1}^\infty$ such that
\begin{enumerate}
\item $\vec{w}_i$ lies on the boundary of $\Kss(X, \sum_{j=1}^kD_j)$ but $\vec{w}_i\notin \{\sum_{j=1}^kx_j=1\}$;
\item $\vec{w}_i$ tends to a point $\vec{w}'$ on the boundary of $\Kss(X, \sum_{j=1}^kD_j)$;
\item  any hyperplane passing through $\vec{w}_i$ and $\vec{w}'$ separates $\Kss(X, \sum_{j=1}^kD_j)$. 
\end{enumerate}
Write $\vec{w}':=(a_1,...,a_k)$ and $\vec{w}_i:=(a_{i1},...,a_{ik})$, then $(X, \sum_{j=1}^k a_{ij}D_j)$ is a K-semistable log Fano pair with $\delta(X, \sum_{j=1}^k a_{ij}D_j)=1$. By Theorem \ref{thm: lxz22} and Theorem \ref{thm: blz22}, there exists a prime divisor $E_i$ over $X$ such that 
$$\beta_{X, \sum_{j=1}^ka_{ij}D_j}(E_i)=A_{X, \sum_{j=1}^ka_{ij}D_j}(E_i)-(1-\sum_{j=1}^ka_{ij})S_X(E_i)=0 ,$$
and $E_i$ induces a special test configuration of $(X, \sum_{j=1}^ka_{ij}D_j)$
$$(\mX_i, \sum_{j=1}^ka_{ij}\mD_{ij})\to \bA^1 $$
such that the central fiber $(\mX_{i,0}, \sum_{j=1}^ka_{ij}\mD_{ij,0})$ is also K-semistable. By Theorem \ref{thm: boundedness} we know that the set $\{(\mX_{i,0}, \sum_{j=1}^k\mD_{ij,0})\}_i$ is log bounded. We first explain that it is enough to show that 
$(\mX_{i,0}, \sum_{j=1}^ka_{j}\mD_{ij,0})$ is a K-semistable log pair (with $\bR$-coefficients) for $i\gg 1$. Suppose 
$(\mX_{i,0}, \sum_{j=1}^ka_{j}\mD_{ij,0})$ is K-semistable for $i\gg 1$. If $\vec{w}'$ does not lie on the hyperplane $\{\sum_{j=1}^kx_j=1\}$, then the central fiber of 
$$(\mX_i, \sum_{j=1}^ka_{j}\mD_{ij})\to \bA^1 $$
is a K-semistable log Fano pair (with $\bR$-coefficients), which implies the following by Lemma \ref{lem: vanishing futaki}:
$$\beta_{X, \sum_{j=1}^ka_{j}D_j}(E_i)=A_{X, \sum_{j=1}^ka_{j}D_j}(E_i)-(1-\sum_{j=1}^ka_{j})S_X(E_i)=0. $$
If $\vec{w}'$ lies on the hyperplane $\{\sum_{j=1}^kx_j=1\}$, then $(\mX_{i,0}, \sum_{j=1}^ka_{j}\mD_{ij,0})$ is an lc log Calabi-Yau pair (with $\bR$-coefficients) since it is K-semistable.  In this case, there exists a finite decomposition
$$(\mX_{i,0}, \sum_{j=1}^ka_j\mD_{ij,0})= \sum_l r_{il}(\mX_{i,0}, \mB_{il,0}:=\sum_{j=1}^kb_{ilj}\mD_{ij,0})$$
which means that $\sum_{j=1}^ka_j\mD_{ij,0}=\sum_l r_{il}\mB_{il,0}$ satisfying
\begin{enumerate}
\item $r_{il}$ are non-negative real numbers with $\sum_l r_{il}=1$;
\item $b_{ilj}$'s are non-negative rational numbers with $\sum_l r_{il}\cdot b_{ilj}=a_j$;
\item $(\mX_{i,0}, \mB_{il,0})$ is an lc log Calabi-Yau pair (with $\bQ$-coefficients) for each $l$.
\end{enumerate}
By \cite[Lemma 2.8]{Zhou23}, $E_i$ is an lc place of $(X, \sum_{j=1}^k b_{ilj}D_j)$ for each $l$, thus $E_i$ is an lc place of $(X, \sum_{j=1}^ka_j D_j)$. This also implies that 
$$\beta_{X, \sum_{j=1}^ka_{j}D_j}(E_i)=A_{X, \sum_{j=1}^ka_{j}D_j}(E_i)-(1-\sum_{j=1}^ka_{j})S_X(E_i)=0. $$
In either case we see that the hyperplane $\beta_{X, \sum_{j=1}^kx_jD_j}(E_i)=0$ passes through both $\vec{w}'$ and $\vec{w}_i$, while it does not separate $\Kss(X, \sum_{j=1}^kD_j)$ by beta criterion. This leads to a contradiction.

We divide the rest of the proof into two parts to show that the log pair $(\mX_{i,0}, \sum_{j=1}^ka_{j}\mD_{ij,0})$ is K-semistable for $i\gg 1$.

\textit{Case 1}. Suppose $\vec{w}'$ is not on the hyperplane $\{\sum_{j=1}^kx_j=1\}.$ We have the following estimation for $\delta(\mX_{i,0}, \sum_{j=1}^ka_{j}\mD_{ij,0})$:
\begin{align*}
&\ \frac{A_{\mX_{i,0}, \sum_{j=1}^ka_{j}\mD_{ij,0}}(F)}{S_{\mX_{i,0}, \sum_{j=1}^ka_{j}\mD_{ij,0}}(F)}\\
=&\ \frac{A_{\mX_{i,0}, \sum_{j=1}^ka_{ij}\mD_{ij,0}}(F)+\sum_{j=1}^k(a_{ij}-a_j)\ord_F(\mD_{ij,0})}{\frac{1-\sum_{j=1}^ka_j}{1-\sum_{j=1}^ka_{ij}}S_{\mX_{i,0}, \sum_{j=1}^ka_{ij}\mD_{ij,0}}(F)}\\
\geq&\ \frac{1-\sum_{j=1}^ka_{ij}}{1-\sum_{j=1}^ka_j}-\frac{\sum_{j=1}^k(d+1)\cdot |a_{ij}-a_j|}{1-\sum_{j=1}^ka_j}.
\end{align*}
It is clear that
$$\frac{1-\sum_{j=1}^ka_{ij}}{1-\sum_{j=1}^ka_j}-\frac{\sum_{j=1}^k(d+1)\cdot |a_{ij}-a_j|}{1-\sum_{j=1}^ka_j}$$
tends to one as $i$ tends to $+\infty$, thus $(\mX_{i,0}, \sum_{j=1}^ka_{j}\mD_{ij,0})$ is K-semistable for $i\gg 1$ by Lemma \ref{lem: real constructibility}.

\

\textit{Case 2}. Suppose $\vec{w}'$ is on the hyperplane $\{\sum_{j=1}^kx_j=1\}$.
Note that $(\mX_{i,0}, \sum_{j=1}^ka_{ij}\mD_{ij,0})$ is a log Fano pair for any $i$. Recall that the set of log pairs $\{(\mX_{i,0}, \sum_{j=1}^k\mD_{ij,0})\}_i$ is log bounded. By Lemma \ref{lem: log bdd}, we see that the log pair $(\mX_{i,0}, \sum_{j=1}^ka_{j}\mD_{ij,0})$ is log canonical for $i\gg 1$ since $\lim_{i\to \infty}a_{ij}=a_j$. This implies that the log pair $(\mX_{i,0}, \sum_{j=1}^ka_{j}\mD_{ij,0})$ is K-semistable for $i\gg 1$.
\end{proof}

\begin{lemma}\label{lem: discrete lemma}
Let $(X, \sum_{j=1}^kD_j)$ be a log pair in $\mE:=\mE(d,k,v,I)$, then the extreme points of $\Kss(X, \sum_{j=1}^kD_j)$ are discrete.
\end{lemma}

\begin{proof}
We assume that $\Kss(X, \sum_{j=1}^kD_j)$ is not contained in the hyperplane $\{\sum_{j=1}^kx_j=1\}$, otherwise the extreme points of $\Kss(X, \sum_{j=1}^kD_j)$ are already discrete by Lemma \ref{lem: lc polytope}. By Lemma \ref{lem: polytope lemma}, $\Kss(X, \sum_{j=1}^kD_j)$ is locally a polytope at every extreme point.
Suppose the extreme points of $\Kss(X, \sum_{j=1}^kD_j)$ are not discrete, then there exists a sequence of points $\{\vec{w}_i\}_{i=1}^\infty$ such that
\begin{enumerate}
\item $\vec{w_i}$ is an extreme point of $\Kss(X, \sum_{j=1}^kD_j)$ for every $i$;
\item $\vec{w}_i$ is not on the hyperplane $\{\sum_{j=1}^kx_j=1\}$ for every $i$;
\item $\vec{w}_i$ tends to a point $\vec{w}$ lying on the boundary of $\Kss(X, \sum_{j=1}^kD_j)$;
\item any hyperplane passing through $\vec{w_i}$ and $\vec{w}$ separates $\Kss(X, \sum_{j=1}^kD_j)$.
\end{enumerate}
Write $\vec{w}_i:=(a_{i1},...,a_{ik})$ and $\vec{w}:=(a_1,...,a_k)$. By Lemma \ref{lem: rationality lemma} we see every $\vec{w}_i$ is rational (while $\vec{w}$ is not necessarily rational). Thus $(X, \sum_{j=1}^ka_{ij}D_j)$ is a K-semistable log Fano pair with $\delta(X, \sum_{j=1}^ka_{ij}D_j)=1$. By Theorem \ref{thm: lxz22} and Theorem \ref{thm: blz22}, there exists a prime divisor $E_i$ over $X$ such that 
$$\beta_{X, \sum_{j=1}^ka_{ij}D_j}(E_i)=A_{X, \sum_{j=1}^ka_{ij}D_j}(E_i)-(1-\sum_{j=1}^ka_{ij})S_X(E_i)=0 ,$$
and $E_i$ induces a special test configuration of $(X, \sum_{j=1}^ka_{ij}D_j)$
$$(\mX_i, \sum_{j=1}^ka_{ij}\mD_{ij})\to \bA^1 $$
such that the central fiber $(\mX_{i,0}, \sum_{j=1}^ka_{ij}\mD_{ij,0})$ is also K-semistable. If $\vec{w}$ is not on the hyperplane $\{\sum_{j=1}^kx_j=1\}$, by the same proof as in case 1 of the proof of Lemma \ref{lem: polytope lemma},  $(\mX_{i,0}, \sum_{j=1}^ka_j\mD_{ij,0})$ is K-semistable for $i\gg 1$ and the hyperplane 
$$\beta_{X, \sum_{j=1}^kx_jD_{j}}(E_i)=0$$ 
passes through $\vec{w}_i$ and $\vec{w}$ while it does not separate $\Kss(X, \sum_{j=1}D_j)$, which is a contradiction. If $\vec{w}$ is on the hyperplane $\{\sum_{j=1}^kx_j=1\}$, then by the same proof as in case 2 of the proof of Lemma \ref{lem: polytope lemma}, $(\mX_{i,0}, \sum_{j=1}^ka_j\mD_{ij,0})$ is an lc log Calabi-Yau pair (with $\bR$-coefficients) for $i\gg 1$. Similarly, there exists a finite decomposition
$$(\mX_{i,0}, \sum_{j=1}^ka_j\mD_{ij,0})= \sum_l r_{il}(\mX_{i,0}, \mB_{il,0}:=\sum_{j=1}^kb_{ilj}\mD_{ij,0})$$
which means that $\sum_{j=1}^ka_j\mD_{ij,0}=\sum_l r_{il}\mB_{il,0}$ satisfying
\begin{enumerate}
\item $r_{il}$ are non-negative real numbers with $\sum_l r_{il}=1$;
\item $b_{ilj}$'s are non-negative rational numbers with $\sum_l r_{il}\cdot b_{ilj}=a_j$;
\item $(\mX_{i,0}, \mB_{il,0})$ is an lc log Calabi-Yau pair (with $\bQ$-coefficients) for each $l$.
\end{enumerate}
By \cite[Lemma 2.8]{Zhou23}, $E_i$ is an lc place of $(X, \sum_{j=1}^k b_{ilj}D_j)$ for each $l$, thus $E_i$ is also an lc place of $(X, \sum_{j=1}^ka_j D_j)$. This implies that the hyperplane $\beta_{X, \sum_{j=1}^kx_jD_{j}}(E_i)=0$ passes through $\vec{w}_i$ and $\vec{w}$ while it does not separate $\Kss(X, \sum_{j=1}D_j)$, a contradiction.
\end{proof}

\begin{theorem}\label{thm: polytope}
Let $(X, \sum_{j=1}^kD_j)$ be a log pair in $\mE:=\mE(d,k,v,I)$, then the K-semistable domain $\Kss(X, \sum_{j=1}^kD_j)$ is a rational polytope.
\end{theorem}

\begin{proof}
Note that the boundary of $\Kss(X, \sum_{j=1}^kD_j)$ is compact.
To show it is a polytope,  it suffices to show that it is locally a polytope at every extreme point and the extreme points of $\Kss(X, \sum_{j=1}^kD_j)$ are discrete. The theorem then follows from Lemma \ref{lem: rationality lemma}, Lemma \ref{lem: polytope lemma}, and Lemma \ref{lem: discrete lemma}.
\end{proof}

\section{Finiteness of polytopes}\label{sec: finiteness}

Fix two positive integers $d$ and $k$, a positive number $v$, and a positive integer $I$. In this section, we will prove that the K-semistable domains of log pairs in $\mE:=\mE(d,k,v,I)$ are finite, which finally leads to the proof of Theorem \ref{thm: main1}.

\begin{lemma}\label{lem: approximation lemma}
Let $(X, \sum_{j=1}^kD_j)$ be a log pair in $\mE$, and $(c_1,...,c_k)$ be a rational vector such that $(X, \sum_{j=1}^kc_jD_j)$ is a log Fano pair. Then there exists a positive number $N$ depending only on the dimension $d$ such that
$\delta(X, \sum_{j=1}^kc_jD_j)$ can be approximated by lc places of $N$-complements of $X$ if $\delta(X, \sum_{j=1}^kc_jD_j)\leq 1$.
\end{lemma}

\begin{proof}
Denote by $\vec{c}:=(c_1,...,c_k)$. We define 
$$\delta_m(\vec{c}):=\delta_m(X,\sum_{j=1}^kc_jD_j):=\inf_E \frac{A_{X, \sum_{j=1}c_jD_j}(E)}{\max_{B_m}(1-\sum_{j=1}^kc_j)\ord_E(B_m)},$$ 
 where $E$ runs through all prime divisors over $X$ and $B_m$ runs through all $m$-basis type divisors of $(X,-K_X)$. This infimum is in fact a minimum, since $m$-basis type divisors form a bounded family and log canonical thresholds take only finitely many values on this family.
 By \cite{BJ20} we have
 $$\lim_m \delta_m(\vec{c})=\delta(X,\sum_{j=1}^kc_jD_j).$$
 
We first assume $\delta(X,\sum_{j=1}^kc_jD_j)<1$, then $\delta_m(\vec{c})<1$ for sufficiently large $m$. In this case one can find an $m$-basis type divisor $B_m$ such that 
$$(X,\sum_{j=1}^kc_jD_j+ \delta_m(\vec{c})(1-\sum_{j=1}^kc_j)B_m)$$ 
is strictly log canonical and admits an lc place $E_m$ over $X$. We choose a general $H\sim_\bQ -K_X$ such that 
$$\left(X,\sum_{j=1}^kc_jD_j+ \delta_m(\vec{c})(1-\sum_{j=1}^kc_j)B_m+a(\vec{c})H\right)$$ 
is an lc log Calabi-Yau pair, where 
$$a(\vec{c}):=1-\sum_{j=1}^kc_j-\delta_m(\vec{c})(1-\sum_{j=1}^kc_j).$$ 
By \cite[Corollary 1.4.3]{BCHM10}, there is an extraction $g_m: Y_m\to X$ which only extracts $E_m$, i.e.
 \begin{align*}
 &\ K_{Y_m}+\sum_{j=1}^kc_j\tilde{D}_j+\delta_m(\vec{c})(1-\sum_{j=1}^kc_j)\tilde{B}_m+a(\vec{c})\tilde{H}+E_m \\
 =\ &\ g_m^*\left(K_X+\sum_{j=1}^kc_jD_j+ \delta_m(\vec{c})(1-\sum_{j=1}^kc_j)B_m+a(\vec{c})H\right),
 \end{align*}
 where $\tilde{\bullet}$ is the birational transformation of $\bullet$. Note that $Y_m$ is of Fano type. We could run $-(K_{Y_m}+E_m)$-MMP to get a mimimal model $f_m: Y_m\dashrightarrow Y_m'$ such that $-(K_{Y_m}+{f_m}_*E_m)$ is nef and $Y_m'$ is also of Fano type. By Theorem \ref{thm: complement}, there exists a positive integer $N$ depending only on $d$ such that $E_m$ is an lc place of some $N$-complement of $X$. Therefore we have
 $$\delta(X,\sum_{j=1}^kc_jD_j)=\lim_m \frac{A_{X,\sum_{j=1}^kc_jD_j}(E_m)}{(1-\sum_{j=1}^kc_j)S_{X}(E_m)},$$
 where $\{E_m\}_m$ is a sequence of lc places of some $N$-complements of $X$.
 
 Next we turn to the case $\delta(X,\sum_{j=1}^kc_jD_j)=1$. By the same proof as \cite[Theorem 1.1]{ZZ21b}, for any rational $0<\epsilon_i\ll1$, one can find an effective $\bQ$-divisor $\Delta_i\sim_\bQ -K_X$ such that 
 $$\delta(X,\sum_{j=1}^kc_jD_j+\epsilon_i\Delta_i)<1.$$ We assume $\{\epsilon_i\}_i$ is a decreasing sequence of rational numbers tending to $0$. For each $\epsilon_i$, by the proof of the previous case, there exists a sequence of prime divisors $\{E_{i,m}\}_m$ which are lc places of some $N$-complements of $X$ such that
 $$\delta(X,\sum_{j=1}^kc_jD_j+\epsilon_i \Delta_i)=\lim_m \frac{A_{X,\sum_{j=1}^kc_jD_j+\epsilon_i\Delta_i}(E_{i,m})}{(1-\sum_{j=1}^kc_j-\epsilon_i)S_{X}(E_{i,m})},$$ 
and $N$ depends only on $d$. Thus we obtain the following 
 $$\delta(X,\sum_{j=1}^kc_jD_j)=\inf_{i,m}\frac{A_{X,\sum_{j=1}^kc_jD_j}(E_{i,m})}{(1-\sum_{j=1}^kc_j)S_{X}(E_{i,m})}. $$
 In both cases, we could find a sequence of lc places of some $N$-complements of $X$ to approximate $\delta(X,\sum_{j=1}^kc_jD_j)$. The proof is finished.
\end{proof}

\begin{lemma}\label{lem: finite lemma}
Let $f: \mX\to B$ be a flat family of Fano varieties of dimension $d$ such that $-K_{\mX/B}$ is a relatively ample $\bQ$-line bundle on $\mX$, and $B$ is a smooth base. Let $\mD_j,\ j=1,...,k$ be effective $\bQ$-divisors on $\mX$ such that every component of $\mD_j$ is flat over $B$ and $\mD_j\sim_{\bQ, B} -K_{\mX/B}$ for every $1\leq j\leq k$. Suppose for every closed point $b\in B$, the fiber $(\mX_b, \sum_{j=1}^k\mD_{j,b})$ is contained in $\mE:=\mE(d,k,v,I)$, and the family $(\mX, \sum_{j=1}^k\mD_j)\to B$ admits a fiberwise log resolution. Then the following  set 
$$\{\Kss(\mX_b, \sum_{j=1}^k\mD_{j,b})\ |\ \textit{$b\in B$}\} $$
is a finite set of K-semistable domains.
\end{lemma}

\begin{proof}

Since the family $(\mX, \sum_{j=1}^k\mD_j)$ admits a fiberwise log resolution, by Lemma \ref{lem: lc polytope}, there exists a rational polytope $P\subset \overline{\Delta^k}$ such that for any rational point $\vec{w}:=(a_1,...,a_k)\in P$, the log pair $(\mX_b, \sum_{j=1}^ka_j\mD_{j,b})$
is log canonical for any closed point $b\in B$. Moreover, if $\vec{w}\in P^\circ$, then $(\mX_b, \sum_{j=1}^ka_j\mD_{j,b})$ is a log Fano pair.

By Lemma \ref{lem: approximation lemma}, there exists a positive integer $N$ depending only on the dimension $d$ such that $\delta(X, \sum_{j=1}^kc_jD_j)$ can be approximated by lc places of $N$-complements of $X$ once we have 
$$\delta(X, \sum_{j=1}^kc_jD_j)\leq 1,$$ 
where $(X, \sum_{j=1}^kD_j)$ is taken from $\mE$ and $(X, \sum_{j=1}^kc_jD_j)$ is a log Fano pair. We emphasize here that the integer number $N$ does not depend on the choice of $(c_1,...,c_k)$.

We may assume that $N$ is sufficiently divisible, then $f_*\mO_\mX(-NK_{\mX/B})$ is a vector bundle over $B$ and we naturally have the following morphism
$$W:=\bP(f_*\mO_\mX(-NK_{\mX/B}))\to B, $$
whose fiber over $b\in B$ is exactly $|-NK_{\mX_b}|$. Let $\mH\subset \mX\times_B W$ be the universal divisor associated to $f_*\mO_\mX(-NK_{\mX/B})$ (this means $\mH_b\subset \mX_b\times |-NK_{\mX_b}|$ is the universal divisor associated to $-NK_{\mX_b}$), and denote by $\mD:=\frac{1}{N}\mH$. Consider the morphism $(\mX, \mD)\to W$. By the lower semicontinuity of log canonical thresholds, the locus
$$Z:=\{w\in W\ |\ \textit{$(\mX_w, \mD_w)$ is log canonical}\} $$
is locally closed in $W$. Then the scheme $Z$ together with the $\bQ$-divisor $\mD_Z$ (which is obtained by pulling back $\mD$ under the morphism $\mX\times_B Z\to \mX\times_B W$) parametrize boundaries of the desired form.

For the morphism 
$$f: (\mX\times_B Z, \mD_Z+\sum_{j=1}^k \mD_j\times_B Z)\to Z,$$ 
we choose a stratification of $Z$, denoted by $Z=\cup_{i=1}^rZ_i$, such that each $Z_i$ is smooth and there is an \'etale cover $Z_i'\to Z_i$ such that 
$$f_i: (\mX\times_B Z, \mD_Z+\sum_{j=1}^k \mD_j\times_B Z)\times _Z Z'_i\to Z'_i$$
admits a fiberwise log resolution.

For any rational point $\vec{w}:=(a_1,...,a_k)\in P^\circ$, we define
$$d_i(\vec{w}):=\inf\Bigg\{\frac{A_{\mX_{z_i}, \sum_{j=1}^ka_j\mD_j\times_B z_i}(v)}{S_{\mX_{z_i}, \sum_{j=1}^ka_j\mD_j\times_B z_i}(v)}\  \Big | \ \textit{$v\in \Val^*_{\mX_{z_i}}$ and $A_{\mX_{z_i}, \mD_{z_i}}(v)=0$}\Bigg\}, $$
where $z_i$ is a closed point on $Z_i'$. By \cite[Proposition 4.2]{BLX22}, the number $d_i(\vec{w})$ is well defined and does not depend on the choice of $z_i\in Z_i'$.  Denote by $\vec{x}:=(x_1,...,x_k)$, then $d_i(\vec{x})$ is a well defined function for $\vec{x}\in P^\circ$. Recall that we use the following notation
$$\tilde{\delta}(X, \Delta):=\min\{\delta(X, \Delta), 1\}$$ 
for a log Fano pair. By Lemma \ref{lem: approximation lemma}, we have
\begin{align*}
 \tilde{\delta}(\mX_b, \sum_{j=1}^k x_j\mD_{j,b})
=\ \big\{\min_i\{1, d_i(\vec{x})\}|\textit{$b$ is contained in the image of $Z_i'\to B$}\big\}. 
\end{align*}
Denote by $p: Z\to B$. We aim to construct a suitable stratification of $B$. Recall that $Z=\cup_{i}^r Z_i$. For each subset $\{i_1,...,i_s\}\subset \{1,2,...,r\}$, we define the following subset of $B$:
$$B_{i_1,...,i_s}:=\Big\{b\in B\ |\ \textit{$b\in p(Z_l)$ if and only if $l\in \{i_1,...,i_s\}$}\Big\}. $$
Clearly we have a finite decomposition of $B$:
$$B=\bigcup_{\{i_1,...,i_s\}\subset\{1,2,...,k\}} B_{i_1,...,i_s}.$$
By our construction, for any $b_1, b_2\in B_{i_1,...,i_s}$, we have
$$\tilde{\delta}(\mX_{b_1}, \sum_{j=1}^k x_j\mD_{j,b_1})=\tilde{\delta}(\mX_{b_2}, \sum_{j=1}^k x_j\mD_{j,b_2})=\min\{1, d_{i_1}(\vec{x}),...,d_{i_s}(\vec{x})\} $$
for any $\vec{x}\in P^\circ$. This implies 
$$\Kss(\mX_{b_1}, \sum_{j=1}^k\mD_{j, b_1})=\Kss(\mX_{b_2}, \sum_{j=1}^k\mD_{j, b_2}).$$ 
Thus there are only finitely many K-semistable domains for fibers of the morphism $(\mX, \sum_{j=1}^k\mD_j)\to B$.
\end{proof}

\begin{theorem}\label{thm: finite}
Fix two positive integers $d$ and $k$, a positive number $v$,  and a positive integer $I$. Then there are only finitely many rational polytopes in $\overline{\Delta^k}$ which may appear as the K-semistable domains of log pairs in $\mE:=\mE(d,k,v,I)$.
\end{theorem}

\begin{proof}
By Theorem \ref{thm: boundedness}, there exists a projective morphism $\mY\to T$ between finite type schemes with a reduced divisor $\mD$ on $\mY$ such that for any $(X, \sum_{j=1}^kD_j)\in \mE$, there exists a closed point $t\in T$ such that 
$$(X, \red(\sum_{j=1}^kD_j))\cong (\mY_t, \mD_t).$$
By the proof of \cite[Theorem 2.21]{XZ20b}, there is a sub-scheme of $T$ which precisely parametrizes the log pairs in $\mE$. Up to a base change, we may assume $(\mY, \mD)\to T$ precisely parametrizes the set $\mE$. Up to a stratification of $T$ and \'etale base changes, we may assume $(\mY, \mD)\to T$ satisfies the conditions in Lemma \ref{lem: finite lemma}. Thus there are only finite polytopes in $\overline{\Delta^k}$ which may appear as K-semistable domains of log pairs in $\mE$. By Theorem \ref{thm: polytope}, these polytopes are rational.
\end{proof}

\begin{proof}[Proof of Theorem \ref{thm: main1}]
The proof is a combination of Theorem \ref{thm: polytope} and Theorem \ref{thm: finite}.
\end{proof}

\begin{proof}[Proof of Corollary \ref{cor: chamber decomposition}]
It is implied by Theorem \ref{thm: main1}.
\end{proof}

\section{Wall crossing for K-moduli}\label{sec: K-moduli}

Fix two positive integers $d$ and $k$, a positive number $v$, a positive integer $I$, and a rational vector $\vec{c}:=(c_1,...,c_k)\in \Delta^k$. Recall that 
$$\Delta^k:=\{(x_1,...,x_k)\ |\ \textit{$x_j\in [0,1)\cap \bQ$ {\rm{and}} $\sum_{j=1}^kx_j<1$}\}.$$
By \cite{Jiang20, Xu20, BLX22}, there exists an Artin stack of finite type, denoted by $\mM^K_{d,k,v,I, \vec{c} }$, parametrizing the log pairs $(X, \sum_{j=1}^kD_j)$ satisfying the following conditions:
\begin{enumerate}
\item $X$ is a Fano variety of dimension $d$ with $(-K_X)^d=v$;
\item $D_j$ is an effective $\bQ$-divisor satisfying $D_j\sim_\bQ -K_X$ for every $j$;
\item $I(K_X+D_j)\sim 0$ for every $j$;
\item $(X, \sum_{j=1}^kc_jD_j)$ is a K-semistable log Fano pair.
\end{enumerate}
By \cite{ABHLX20, BHLLX21, CP21, XZ20b, LXZ22}, the stack $\mM^K_{d,k,v,I, \vec{c}}$ descends to a projective scheme $M^K_{d,k,v,I,\vec{c}}$ as a good moduli space which parametrizes K-polystable objects. We call $\mM^K_{d,k,v,I,\vec{c}}$ (resp. $M^K_{d,k,v,I,\vec{c}}$) a K-moduli stack (resp. K-moduli space). By the finite chamber decomposition established in Theorem \ref{thm: main1} and Corollary \ref{cor: chamber decomposition}, we have the following wall crossing for K-moduli.

\begin{theorem}
Fix two positive integers $d$ and $k$, a positive number $v$, and a positive integer $I$. There exists a finite chamber decomposition $\overline{\Delta^k}=\cup_{i=1}^r P_i$, where $P_i$ are rational polytopes and $P_i^\circ\cap P_j^\circ=\emptyset$ for $i\ne j$, such that  
\begin{enumerate}
\item For any face $F$ of any chamber $P_i$, $\mM^K_{d,k,v,I,\vec{c}}$ {\rm{(resp. $M^K_{d,k,v,I,\vec{c}}$)}} does not change as $\vec{c}$ varies in $F^\circ(\bQ)$; 
\item Let $P_{i}$ and $P_{j}$ be two different polytopes in the decomposition that share the same face $F_{ij}$. Suppose $\vec{w}_1\in P_i^\circ(\bQ)$ and $\vec{w}_2\in P_j^\circ(\bQ)$ satsify that the segment connecting them
intersects $F_{ij}$ at a point $\vec{w}$, then $\vec{w}$ is a rational point and we have the following diagram for any rational numbers $0<t_1, t_2< 1$:
\begin{center}
	\begin{tikzcd}[column sep = 1em, row sep = 2em]
	 \mM_{d,k,v,I, t_1\vec{w}+(1-t_1)\vec{w}_1}^K \arrow[d,"",swap] \arrow[rr,""]&& \mM_{d,k,v,I, \vec{w}}^K\arrow[d,"",swap] &&\mM_{d,k,v,I, t_2\vec{w}+(1-t_2)\vec{w}_2}^K\arrow[d,""]\arrow[ll,""]\\
	 M_{d,k,v,I, t_1\vec{w}+(1-t_1)\vec{w}_1}^K\arrow[rr,""]&& M_{d,k,v,I, \vec{w}}^K&&M_{d,k,v,I, t_2\vec{w}+(1-t_2)\vec{w}_2}^K\arrow[ll,"", swap].
	 	 	 	\end{tikzcd}
\end{center}
Moreover, we have
\begin{enumerate}
\item $\mM_{d,k,v,I, t_l\vec{w}+(1-t_l)\vec{w}_l}^K$
{\rm{(resp. $M_{d,k,v,I, t_l\vec{w}+(1-t_l)\vec{w}_l}^K$)}} does not change as $t_l$ varies in $[0,1)\cap \bQ$, where $l=1,2$;
\item $ \mM_{d,k,v,I, t_1\vec{w}+(1-t_1)\vec{w}_1}^K\longrightarrow \mM_{d,k,v,I, \vec{w}}^K\longleftarrow \mM_{d,k,v,I, t_2\vec{w}+(1-t_2)\vec{w}_2}^K$ are open embeddings for $t_1, t_2\in [0,1)\cap \bQ$;
\item $ M_{d,k,v,I, t_1\vec{w}+(1-t_1)\vec{w}_1}^K\longrightarrow M_{d,k,v,I, \vec{w}}^K\longleftarrow M_{d,k,v,I, t_2\vec{w}+(1-t_2)\vec{w}_2}^K$ are projective morphisms for $t_1,t_2\in [0,1)\cap \bQ$.
\end{enumerate}
\end{enumerate}
\end{theorem} 

\begin{proof}
The first statement follows from Theorem \ref{thm: main1} and Corollary \ref{cor: chamber decomposition}. For the second statement, the rationality of $\vec{w}$ follows from Lemma \ref{lem: dense lemma}, (2.a) follows from (1), (2.b) follows from \cite{BLX22, Xu20}, and (2.c) follows from \cite{XZ20b, LXZ22}.
\end{proof}

\begin{remark}
For general type varieties, the wall crossing theory for KSBA-moduli with multiple boundaries is studied in \cite{ABIP21}.
\end{remark}

\begin{remark}
There are two further directions worth consideration on this topic. The first one is to work out explicit examples on the chamber decomposition of K-semistable domains. For example, compute all the K-semistable domains that may appear for the log pair $(\bP^n, Q_1+Q_2)$ as $Q_i, i=1,2,$ vary in the linear system $|\mO_{\bP^n}(2)|$ (note that we have confirmed the K-semistable domain when $(\bP^n, Q_1+Q_2)$ is log smooth, e.g \cite{Zhou24, LZ23}). The second is to clarify the wall crossing theory (in particular the nature of the K-semistable domain) in non-linear setting. Here \textit{non-linear} means that the boundary divisors of the pair are not necessarily proportional to the anti-canonical bundle.
\end{remark}

\bibliography{reference.bib}
\end{document}